\documentclass[12pt,a4paper]{article}
\usepackage{amsmath,amsthm,amsfonts,amssymb,latexsym,xy}
\xyoption{all} 
\usepackage{enumerate}
\usepackage{a4wide}
\newtheorem{theorem}{Theorem}[section]
\newtheorem{definition}[theorem]{Definition}
\newtheorem{lemma}[theorem]{Lemma}
\newtheorem{remark}[theorem]{Remark}
\newtheorem{proposition}[theorem]{Proposition}
\newtheorem{cor}[theorem]{Corollary}
\newtheorem{example}[theorem]{Example}
\newcommand{\supp}{\text{supp }}
\newcommand{\im}{\text{Im }}
\newcommand{\Y}{{\mathfrak Y}}
\newcommand{\X}{{\mathfrak X}}
\newcommand{\tr}{\text{tr}}
\renewcommand{\S}{{\mathfrak S}}
\renewcommand{\H}{{\mathcal H}}

\begin{document}
\date{}

\title{Operator synthesis II. Individual synthesis and linear operator 
equations}
\author{{\it Victor Shulman} at Vologda and 
{\it Lyudmila  Turowska} at G\"oteborg}

\maketitle{}
\begin{abstract}
The second part of our work on operator synthesis deals with individual
operator synthesis
of elements in some tensor products, in particular in Varopoulos algebras, and
its connection with linear operator equations. Using a developed technique
of ``approximate inverse intertwining'' we obtain some generalizations
of the Fuglede and the Fuglede-Weiss theorems and solve some problems posed in
\cite{o,weiss1,weiss2}. Additionally, we give some 
applications to spectral synthesis in Varopoulos algebras and to partial
differential equations.
\end{abstract}

\section{Introduction}
This work is a sequel of \cite{sht}  where the problems of operator synthesis 
were treated ``globally'' for lattices of subspaces, bilattices, or, in 
coordinate setting, for subsets of direct products of measure spaces. Here we
consider operator-synthetic properties of elements of some tensor products, 
first of all of the Varopoulos algebras $V(X,Y)=C(X)\hat\otimes C(Y)$. The 
topic  is deeply connected to the theory  of linear operator equations and,
more generally, to the spectral theory of multiplication operators
in the space of bounded operators and in symmetrically normed ideals
of operators. We obtain some extensions of the Fuglede and Fuglede-Weiss theorems,
answer several questions posed in \cite{o,weiss1,weiss2}, give applications to
spectral synthesis in Varopoulos algebras and (somewhat unexpectedly) to
partial differential equations.

Let us describe the results of the paper in more detail. In Section~2 we 
consider some pseudo-topologies and functional spaces on direct products of 
measure spaces. Basic definitions and results from \cite{arv} and \cite{sht}
 related to operator synthesis  for subsets in a direct product  $X\times Y$
 are recalled. It is proved that a subset with a scattered family of $X$-sections 
is equivalent to a countable union of rectangles. 
A consequence which is used later on is that
the set of all solutions $(x,y)$ of an equation of the form 
$\sum_{i=1}^n a_i(x)b_i(y)=0$ is a union of countable family of rectangles 
and a  set of measure null. A special case of this result was established in 
 \cite[Proposition~12]{weiss1}. 

Section~3 deals with a kind of spectral synthesis in commutative Banach 
algebras -- the synthesis with respect to Banach modules. The real distinction 
of this
theory from the classical one is that given a module we get a special class of 
ideals, the annihilators of subsets in the module, and  work  with them only.
 Here our aim is to compare the conditions for an element to be synthetic with
respect to a module and to admit spectral
 synthesis in the algebra. We also relate these conditions to spectra and 
spectral subspaces of the corresponding multiplication operators.

In Section~4 the approach is reduced to the case of  operator modules over
Varopoulos algebras. Let $\mu$, $\nu$ be regular measures on compacts $X$, $Y$
 and $H_1$, $H_2$  
 the corresponding $L_2$-spaces. Then the space $B(H_1, H_2)$ of all bounded 
operators from $H_1$ to $H_2$ becomes a $V(X,Y)$-module with respect to the 
action $(f\otimes g)\cdot T = M_gTM_f$, where $M_f$, $M_g$ are the multiplication 
operators. It is proved that  $F\in V(X,Y)$ admits spectral synthesis iff 
it is synthetic with respect to all modules of this kind. This allows to 
obtain results on spectral synthesis in an  operator-theoretical way.  
The following auxiliary statement (Corollary~\ref{co}) 
appears to be useful: a function $F\in V(X,Y)$ is synthetic with respect to  
$B(H_1,H_2)$
 iff the space of all solutions of the equation $F\cdot X=0$ is reflexive (in 
the sense of \cite{lsh}) 
and iff the $0$-spectral subspace of the operator of multiplication by $F$
coincides with its kernel. 

 The topics of Section~5 are linear operator equations of general type and 
modules over weak$^*$-Haagerup tensor products of $L_{\infty}$-algebras. 
Some estimates for the action of a linear multiplication operator with normal 
coefficients on its $0$-spectral subspace are obtained. 
The extension of the approach allows to relate the topic with global operator 
synthesis.

In Section~6 we develop a general technique which relates solutions
of the ''same'' linear equations in different linear topological spaces.
More strictly speaking we are given  two operators, $S$ and $T$, acting in spaces $\X$, $\Y$, 
and a linear injection, $\Phi: \X \to \Y$, that intertwines them: $T\Phi = \Phi S$.
Our  main tool then is the ``approximate inverse intertwining'' (AII), 
that 
is a net $\{F_{\alpha}\}$ of maps from $\Y$ to $\X$ satisfying the 
conditions that $F_{\alpha}\Phi\to 1_{\X}$, $\Phi F_{\alpha}\to 1_{\Y}$ and
$F_{\alpha}T-SF_{\alpha}\to 0_{\X}$ in the topology of simple convergence.
It appears to be possible to obtain some non-completely trivial results on 
inclusions of images or norm-inequalities in such a general abstract scheme.
  
In applications of the AII-technique to linear operator  equations, 
the spaces $\X$, $\Y$ are symmetrically normed ideals of the algebra $B(H)$ 
(actually the case $\Y = B(H)$ is the most important) and $S$, $T$ are the 
restrictions of a multiplication operator $\Delta$ to $\X$, $\Y$. 
The conditions 
under which an AII exists are considered in Section~7. They are 
close in spirit to  Voiculescu's conditions of quasidiagonality modulo a 
symmetrically normed ideal, but  formally are more weak: instead of the 
condition $||[A,P_n]|| \to 0$ we need only the boundedness of the norms 
$||[A,P_n]||$ (semidiagonality). Note that 
for the usual operator norm the semidiagonality holds automatically while 
quasidiagonality is an intriguing property which was explored in a great  
number of  publications. We discuss examples of  $\S_p$-quasidiagonal 
families; the most simple ones are families of  weighted shifts ($p = 1$) 
and families of commuting normal operators with thin joint spectra.

In Section~8 the applications of  AII's to the problem of triviality of  
the trace of a commutator  and to some related problems are gathered.  
In \cite{weiss} Weiss proved that if a commutator $[A,X]$ of a normal operator 
$A$ and a Hilbert-Schmidt operator $X$ belongs to $\S_1$ then $\tr([A,X])=0$. 
In Proposition~\ref{aiipr5} we extend this result as follows: if a family, 
$\{A_k\}_{k=1}^n$ of operators is $\S_{p/(p-1)}$-semidiagonal  and if a sum 
$\sum_{k=1}^n[A_k,X_k]$
 belongs to $\S_1$ then   $\tr(\sum_{k=1}^n[A_k,X_k])=0$.  We answer also 
 some questions posed in \cite{weiss} and \cite{o} which  are 
formulated in purely function-theoretical terms but in their essence are about  the trace of (sums of) commutators.

 The famous Fuglede Theorem can be formulated as the equality 
$\ker\Delta =\ker\tilde\Delta~$, where $\Delta(X)  = AX-XA$, 
$\tilde\Delta(X)=A^*X-XA^*$ and $A$ is a normal operator. Weiss \cite{weiss} 
strengthen the result to $||\Delta(X)||_2=||\tilde\Delta(X)||_2$. 
Weiss also proposed to consider the case when $\Delta$ is a more general 
multiplication operator $\Delta(X) = \sum_{k\in K} B_kXA_k$ with commuting normal 
coefficients and $\tilde\Delta(X)=\sum_{k\in K} B_k^*XA_k^*$.  
He proved the 
equality $||\Delta(X)||_2=||\tilde\Delta(X)||_2$ in the case when $K$ is 
finite and  both parts of 
the equality  are finite (that is $\Delta(X)$ and $\tilde\Delta(X)$  
belong to $\S_2$).  In general,  these restrictions can not be dropped 
(\cite{zap}).
We show that if the Hausdorff dimension of the joint spectrum of the 
family $\{A_k\}$
 does not exceed $2$ the equality holds without the restrictions. 
 We discuss also a ``non-commutative version'' of the Fuglede Theorem: 
$\ker\tilde\Delta\Delta = \ker\Delta$, where $\Delta$ is arbitrary 
(the coefficients are not supposed to be normal or commuting). 
It is proved in Theorem~\ref{aiicor6} that  the equality 
holds  if $\{A_k\}_{k\in K}$ is 
$1$-semidiagonal. 

 In Theorem~\ref{aiipr} we show the inequality $||\Delta(X)||_2 \geq 
||\tilde\Delta(X)||_2$ for  $\Delta(X) = AX-XB$, provided that  $A$ and $B^*$ are 
hyponormal operators of finite multiplicity.

 Section~10 is devoted to multiplication operators with normal finite 
families of coefficients. It is proved that the ascent of such an operator 
does not exceed $d/2$ where $d$ is the Hausdorff dimension of the
 joint spectrum of  the left coefficient family. This result is applied in 
Section~11 to  the  evaluation of the number on which  the chain of closed 
ideals generated by the powers of an element $F$ of  a Varopoulos algebra is stabilized. 
In particular, it is proved that if $F(x,y)=\sum_{k=1}^nf_k(x)g_k(y)\in 
V(X,Y)$,  
$\dim X \leq 2$  and the functions $f_k$  are Lipschitsian then $F$ admits 
 spectral synthesis. 

 Our last application is  to partial differential equations
 with constant 
coefficients.  Corollary~\ref{diffeq} states that the space of all 
bounded solutions of the equation 
$p(i\frac{\partial}{\partial x_1},i\frac{\partial}{\partial x_2})u=0$ depends 
only on the variety of zeros of the polynomial $p$.

The authors are grateful to Gary Weiss for stimulating results, questions 
and discussions. 

The work was partially written when the first author was visiting 
Chalmers University of Technology in G\"oteborg, Sweden.
The research  was partially supported by a grant from the Swedish Royal 
Academy of Sciences as a part of the program of cooperation with former
Soviet Union. 

\section{Pseudo-topologies and functional spaces on direct products of
 measure spaces}\label{pseudo}

Let $(X,\mu)$, $(Y,\nu)$ be  standard measure spaces with
finite measures,
$m=\mu\times\nu$ the product measure on $X\times Y$. In this section we 
recall some definitions and results from \cite{arv,eks,sht} and obtain  a
few others auxiliary results.

A rectangle in $X\times Y$ is a measurable subset of the form $A\times B$,
where $A\subset X$, $B\subset Y$.
 
%Let $M(X)$ denote the set of finite Borel measures on a compact $X$.
%For $(\mu,\nu)\in M(X)\times M(Y)$, set $H_1=L_2(X,\mu)$, $H_2=L_2(Y,\nu)$.
A subset $E\subset X\times Y$ is called {\it marginally null}  (with respect 
to $\mu\times\nu$) if $E\subset (X_1\times Y)\cup(X\times Y_1)$ and
$\mu(X_1)=\nu(Y_1)=0$.
Two subsets $E_1$, $E_2$ are marginally equivalent ($E_1\sim ^M E_2$ or 
simply $E_1\cong E_2$) if
their symmetric difference is marginally null. Furthermore,
$E_1\subset^M E_2$ means that $E_1\setminus E_2$ is marginally null, a 
property holds marginally almost everywhere if it holds everywhere apart
of a marginally null set, and so on.

A subset  $E$ is called {\it pseudo-open} (more strictly, 
$\omega$-pseudo-open) if it is marginally equivalent to a countable union 
of measurable rectangles.
The complements of pseudo-open sets
are {\it pseudo-closed} sets.  

It is easy to see that the family of all pseudo-open sets defines a 
pseudo-topology on $X\times Y$: it is stable under finite intersections
and countable unions. This pseudo-topology is denoted by $\omega$.

A complex-valued function $f$ on $X\times Y$ is pseudo-continuous if
$f$-preimages of open sets are pseudo-open. It is known (\cite{eks}) that
pseudo-continuous functions form a functional algebra on $X\times Y$. 
In particular, all functions of finite length $f(x,y)=\sum_{i=1}^na_i(x)b_i(y)$
(with measurable $a_i$, $b_i$) are pseudo-continuous.
  
Set  $\Gamma(X,Y)=L_2(X,\mu)\hat\otimes L_2(Y,\nu)$, where $\hat\otimes$ 
denotes the projective tensor product. Clearly, every 
$\Psi\in\Gamma(X,Y)$ can be identified with a function
$\Psi:X\times Y\to{\mathbb C}$ which
admits a representation
\begin{equation}\label{eq1}
\Psi(x,y)=\sum_{n=1}^{\infty}f_n(x)g_n(y)
\end{equation}
where $f_n\in L_2(X,\mu)$, $g_n\in L_2(Y,\nu)$   and 
$\sum_{n=1}^{\infty}||f_n||_{L_2}\cdot||g_n||_{L_2}<\infty$.
 Such a representation
defines a function
marginally almost everywhere (m.a.e.), so two functions in  $\Gamma(X,Y)$ which coincides 
m.a.e. are identified.
$L_2(X,\mu)\hat\otimes L_2(Y,\nu)$-norm of $\Psi$ is
$$||\Psi||_{\Gamma}=\inf\sum_{n=1}^{\infty}||f_n||_{L_2}\cdot
||g_n||_{L_2},$$ 
where the infimum is taken over all  sequences $f_n$, $g_n$ for which
(\ref{eq1}) holds m.a.e. 

We consider also the space $V^{\infty}(X,Y)$ of all (marginal equivalence 
classes of) functions $\Psi(x,y)$ that can be written in the form
(\ref{eq1}) with $f_n\in L^{\infty}(X,\mu)$, $g_n\in L^{\infty}(Y,\nu)$ and
$$\sum_{n=1}^{\infty} |f_n(x)|^2\leq C,\quad x\in X,\quad \sum_{n=1}^{\infty}|g_n(y)|^2\leq C,\quad 
y\in Y.$$
The least possible $C$ here is the norm of $\Psi$ in $V^{\infty}(X,Y)$.
In tensor notations 
$V^{\infty}(X,Y)=L^{\infty}(X,\mu)\hat\otimes^{w^*h} L^{\infty}(Y,\nu)$,
the weak$^*$-Haagerup tensor product  
(\cite{BS}, see also \cite{spronk}, where these are called measurable Schur multipliers). Since  measures $\mu$, $\nu$ are finite, $V^{\infty}(X,Y)\subset 
\Gamma(X,Y)$.
\begin{lemma}\label{EKS}\cite{eks}
All functions $\Psi\in\Gamma(X,Y)$ are $\omega$-pseudo-continuous.
\end{lemma}

Now we discuss the null sets of (families of) functions in $\Gamma(X,Y)$.
We say that $F\in\Gamma(X,Y)$ vanishes on $E\subset X\times Y$ if
$F\chi_{E}=0$ m.a.e., $\chi_E$ being the characteristic function of $E$. 
For ${\mathcal F}\subset\Gamma(X,Y)$, the {\it null set}, 
$\text{null }{\mathcal F}$,  is defined 
to be the largest, up to marginal equivalence, pseudo-closed set such that
each function $F\in{\mathcal F}$ vanishes on it. 
If $E$ is a pseudo-closed subset of $X\times Y$,
let
$$\Phi(E)=\{F\in\Gamma(X,Y)\mid F\text{ vanishes on }E\},$$
$$\Phi_0(E)=\overline{\{F\in\Gamma(X,Y)\mid F\text{ vanishes on a 
nbhd of }E\}},$$
where by a neighborhood  we mean a pseudo-open set containing $E$ and the 
closure is taken in $\Gamma(X,Y)$. 
By \cite{sht}[Theorem~2.1], $\Phi_0(E)$ and $\Phi(E)$ are the smallest and the
largest invariant (with respect to the multiplication by functions
$f\in L_{\infty}(X,\mu)$ and $g\in L_2(Y,\nu)$) closed subspaces of
$\Gamma(X,Y)$ whose null
set is $E$.

We will also need  another pseudo-topology on $X\times Y$. Let us say that
a subset $E\subset X\times Y$ is $\tau$-pseudo-open if it is a union of an
$m$-null set and a countable family of rectangles. It is not difficult to 
check that the class of all such sets is stable under finite intersections 
and countable unions. Clearly, the pseudo-topology $\tau$ is stronger than 
$\omega$. In particular, all functions of finite length and functions in
$\Gamma(X,Y)$ are $\tau$-pseudo-continuous.

Our next aim is to obtain some sufficient condition for a set to be
$\tau$-pseudo-open.

For $U$, $V\subset X$ let us write $U\subset^{\mu}V$  if $\mu(U\setminus V)=0$.
If $U\subset^{\mu}V$ and $V\subset^{\mu}U$ we say that $U$ and $V$ are 
$\mu$-equivalent and write $U\sim^{\mu} V$.

A family ${\mathcal F}$ of measurable subsets of $X$ is called $\mu$-scattered 
if $any$ decreasing sequence $U_1\supset U_2\supset\ldots$ of its members 
$\mu$-stabilizes (that is $U_n\sim^{\mu} U_{n+1}\sim^{\mu}U_{n+2}\sim\ldots$
for some $n$). 

Let now $E$ be a subset of $X\times Y$. An $X$-section of $E$ is a subset of the form
$$E^y=\{x\in X\mid (x,y)\in E\}.$$
$E$ is called {\it $X$-scattered} if the family of all finite intersections of its $X$-sections is $\mu$-scattered. The following result gives us important examples of $X$-scattered sets.
\begin{proposition}
Let $h(x,y)$ be a complex-valued function on $X\times Y$ that has a finite 
length, that is $h(x,y)=\sum_{i=1}^Na_i(x)b_i(y)$, where $a_i$ and $b_i$ are
measurable. Then
$$E=\{(x,y)\mid h(x,y)=0\}$$ is an $X$-scattered set.
\end{proposition}
\begin{proof}
Let $\vec{a}:X\to{\mathbb C}^N$,  $\vec{b}:X\to{\mathbb C}^N$ be defined by
$\vec{a}(x)=\{a_i(x)\}_{i=1}^N$, $\vec{b}(y)=\{b_i(y)\}_{i=1}^N$.
For any $y\in Y$, the $X$-section $E^y$ is the preimage with respect
to $\vec{a}$ of the hyperplane 
$$\{\vec{z}\in{\mathbb C}^N\mid\sum_{i=1}^Nz_ib_i(y)=0\}.$$
Since intersections of hyperplanes must stabilize, the family of their preimages is scattered.
\end{proof}

It is easy to see that any set, which is $m$-equivalent to a finite union of 
rectangles, is  $X$-scattered.
\begin{theorem}\label{scat}
Any $X$-scattered set is $m$-equivalent to a countable union of rectangles.
\end{theorem}
\begin{proof}
Let $E\subset X\times Y$ be an  $X$-scattered set. Denote by  ${\mathcal U}$ the set of 
all countable unions of rectangles 
 and set
$$\underline{m}(E)=\sup\{m(U)\mid U\in{\mathcal U}, U\subset^m E\}.$$
Choosing $U_n\in{\mathcal U}$ with $U_n\subset^m E$ and 
$\displaystyle m(U_n)>\underline{m}(E)-\frac{1}{n}$, we set $U=\cup_{n=1}^{\infty}U_n$.
Then $m(U)=\underline{m}(E)$, $U\in{\mathcal U}$, $U\subset^m E$.
Hence the set $S=E\setminus U$ has the property
that $m(S\cap\Pi)=0$ for any rectangle $\Pi\subset^m E$.
It remains to show that $m(S)=0$.

We define a measure $\nu$ on $X$ by $\nu(A)=m(S\cap(A\times Y))$.
Clearly, $\nu\ll\mu$ whence by the Jordan Theorem $X=X_0\cup X_1$, 
$X_0\cap X_1=\emptyset$,
$\nu(X_0)=0$ and $\nu\sim\mu$ on $X_1$. If $\mu(X_1)=0$ then $m(S)=0$ and we are done.

Assume $\mu(X_1)\ne 0$ and 
let $E_1=E\cap(X_1\times Y)$. Clearly  $E_1$ is $X$-scattered. 
Let $A_1\supset A_2\supset\ldots\supset A_n$ be a maximal chain in the family of
(the equivalence classes of) finite intersections of $X$-sections of $E_1$, 
and let $A$ be its final non-zero element ($A=A_n$ if $\mu(A_n)\ne 0$, 
$A=A_{n-1}$ otherwise). Then any $X$-section of
$E_1$ either $\mu$-contains $A$ or is $\mu$-disjoint with $A$.
Set 
$$K=\{y\in Y\mid A\subset^{\mu}E_1^y\}.$$
Then all $X$-sections of the set $(A\times(Y\setminus K))\cap E_1$
are $\mu$-null. Therefore $m((A\times(Y\setminus K))\cap E_1)=0$,
 $m((A\times(Y\setminus K))\cap S)=0$ and 
$m((A\times K)\cap S)=m((A\times Y)\cap S)=\nu(A)\ne 0$ 
(because $\mu(A)\ne 0$).

On the other hand, $(A\times K)\setminus E_1$ has $\mu$-null sections 
whence 
$m((A\times K)\setminus E_1)=0$. Thus the rectangle $\Pi=A\times K$ is 
$m$-contained in $E_1$ and has non-trivial intersection with $S$, a contradiction.
\end{proof}

\begin{cor}
A pseudo-closed $X$-scattered set is $\tau$-pseudo-open.
\end{cor}
\begin{proof}
Let $E$ be a pseudo-closed $X$-scattered set. If $\Pi\subset^m E$, where $\Pi$ is a rectangle, then $\Pi\subset^ME$. Indeed, $m(\Pi\cap\Pi')=0$ for each rectangle $\Pi'\subset E^c$ so
$\Pi\cap\Pi'\sim^M 0$.

It follows that $\Pi$ can be changed by a sub-rectangle $\tilde\Pi$ of the same 
measure such that $\tilde\Pi\subset E$. This clearly implies our statement: since
$E\sim^m\cup_{j=1}^\infty\Pi_j$ we have $E\sim^m\cup_{j=1}^{\infty}\tilde\Pi_j\subset E$.
\end{proof}
\begin{cor}\label{scatcor}
The set of zeros of a function of finite length is $\tau$-pseudo-open. 
\end{cor}
Now it is easy to deduce a more general result.
\begin{cor}
Let $h_j$, $1\leq j\leq n$, be real-valued functions of finite length on
$X\times Y$. The  set
$$E=\{(x,y)\mid h_j(x,y)\leq 0,1\leq j\leq n\}$$
is $\tau$-pseudo-open. 
\end{cor}
\begin{proof}
Clearly $E=\cap_{j=1}^nE_j$, where $E_j=\{(x,y)\mid h_j(x,y)\leq 0\}$,
and $E_j=E_j'\cup E_j''$, where 
$$E_j'=\{(x,y)\mid h_j(x,y)= 0\},\quad E_j''=\{(x,y)\mid h_j(x,y)< 0\}.$$
All $E_j'$ are $\tau$-pseudo-open by Corollary~\ref{scatcor}. All $E_j''$ are
$\omega$-pseudo-open by Lemma~\ref{EKS} and hence $\tau$-pseudo-open.
\end{proof}

For us the space $\Gamma(X,Y)$ is important because it is predual to the space 
of bounded operators, $B(H_1,H_2)$,  from 
$H_1=L_2(X,\mu)$ to $H_2=L_2(Y,\nu)$.
The duality is given by 
 $$\langle T,\Psi\rangle=\sum_{n=1}^{\infty}(Tf_n,\bar{g}_n),$$ 
with $T\in B(H_1,H_2)$ and 
$\Psi(x,y)=\sum_{n=1}^{\infty}f_n(x)g_n(y)$.

Let $P_U$ and $Q_V$ denote the multiplication operators by 
the characteristic functions of $U\subset X$ and $V\subset Y$.
We say  that $T\in B(H_1,H_2)$ is {\it supported} in 
$E\subset X\times Y$ (or $E$ supports $T$) if $Q_VTP_U=0$ for each Borel sets 
$U\subset X$, 
$V\subset Y$ such that $(U\times V)\cap E=\emptyset$.  
 Then there exists the smallest (up to a marginally null set)
pseudo-closed set, $\supp T$, which supports $T$.
More generally, for any subset ${\mathfrak M}\subset B(H_1,H_2)$ there is the smallest 
pseudo-closed set $\supp{\mathfrak M}$, which supports all operators in ${\mathfrak M}$.
In the  seminal paper \cite{arv} Arveson  defined a support in a similar way 
but using closed sets instead of pseudo-closed (in his setting $X$, $Y$ are
topological spaces). This closed support, $\text{supp}_A T$ can be strictly larger 
than $\supp T$.

For any pseudo-closed set $E\subset X\times Y$ the set, ${\mathfrak M}_{max}(E)$,
of all operators $T$, supported in $E$, has support $E$ and is the largest set 
with this property. It is easy to check that ${\mathfrak M}_{max}(E)$ is a
${\mathcal D}_1\times{\mathcal D}_2$-bimodule, where 
${\mathcal D}_1$, ${\mathcal D}_2$  are the algebra of multiplications 
by functions in  $L_{\infty}(X,\mu)$ and $L_{\infty}(Y,\nu)$ respectively.
There is also the smallest bimodule 
${\mathfrak M}_{min}(E)\subset B(H_1,H_2)$ with
support equal to  $E$ and, moreover,
$${\mathfrak M}_{max}(E)=\Phi_0(E)^{\perp}, \quad 
{\mathfrak M}_{min}(E)=\Phi(E)^{\perp}.$$
(see \cite{sht}).
We say that a pseudo-closed set $E\subset X\times Y$ is 
{\it operator synthetic} (or
$\mu\times\nu$-synthetic) if  the following equivalent conditions hold:
\begin{itemize}
\item $\Phi(E)=\Phi_0(E)$.
\item ${\mathfrak M}_{max}(E)={\mathfrak M}_{min}(E)$.
\item $\langle T,F\rangle=0$ for any $T\in B(H_1,H_2)$ and
$F\in\Gamma(X,Y)$ with $\supp T\subset E\subset\text{null }F$.
\end{itemize}
In further sections we will use also the following  characterization
of  ${\mathfrak M}_{min}(E)$ 
by pairs of projections. 
Identifying projections $P\in B(l_2)\bar\otimes{\mathcal D}_1$
and $Q\in B(l_2)\bar\otimes{\mathcal D}_2$  with projection-valued functions
$P(x):X\to B(l_2)$ and $Q(y):Y\to B(l_2)$ we say that a pair $(P,Q)$ is
an $E$-pair if $P(x)Q(y)$ vanishes on $E$ m.a.e.  
Then by \cite{sht}[Corollary~4.4]
$${\mathfrak M}_{min}(E)=\{T\in B(H_1,H_2)\mid Q(1\otimes T)P=0\text{ for any 
$E$-pair }(P,Q)\}.$$

\section{Synthesis with respect modules over Banach algebras}
Let $A$ be a semisimple, regular, commutative Banach algebra with unit and 
let $X_{A}$ be its spectrum. For any $a\in A$ we shall denote by 
$\hat a$ its Gelfand transform and set
$$null(a)=\{\chi\in X_A\mid\hat a(\chi)=0\}.$$
More generally, for any subset $B \subset A$, we define 
$null(B)$ as $\displaystyle\cap_{a\in B}null(a)$. 

To any closed subset 
$E\subset X_A$ there correspond ideals  
\begin{eqnarray*}
&I(E)= \{r \in A\mid \hat{r}^{-1}(0) \text{ contains } E\},\\
&J_0(E)=\{r\in A\mid \hat{r}^{-1}(0) \text{ contains a nbhd of }E\}
 \text{ and }
J(E)=\overline{J_0(E)}.
\end{eqnarray*}
It is known that $null(J(E)) = null(I(E)) = E$ and 
$J(E) \subset K \subset I(E)$, for any closed ideal $K$ with $null(K) = E$.

For $a\in A$ we define
\begin{eqnarray*}
&I_a= I(null(a)),\\
&J_a^0= J_0(null(a)), \ \mbox{and}\ J_a= J(null(a)).
\end{eqnarray*}
One says that $a\in A$ {\it admits spectral synthesis} if $a\in J_a$. 

Let $M$ be a Banach $A$-module. For any $x\in M$ set
\begin{eqnarray*}
&ann(x)=\{a\in A\mid a\cdot x=0\},\\
&Supp(x)= null(ann(x)).
\end{eqnarray*}
Then $ann(x)$ is a closed ideal and $Supp(x)$ is a closed subset in 
$X_A$.  

In a similar way one defines $ann(N)$ and $Supp(N)$ for arbitrary
subset $N\subset M$.
\begin{definition}\rm
We say that an element $a\in A$ {\it admits synthesis} (or {\it IS synthetic}) 
with respect to an $A$-module $M$
if $a\cdot x=0$ for any $x\in M$ such that $Supp (x)\subset null(a)$. 
\end{definition}

\begin{example}\rm
$A$ is a module over itself with the action defined by 
$a\cdot x=ax$ for any $a,x\in A$. Each $a\in A$ admits synthesis with respect
 to $A$ as  
$A$-module. In fact, assuming $ax\ne 0$ for some $a,x\in A$ with 
$Supp(x)\subset null(a)$, one can find  $\chi\in X_A$ such that 
$\hat{ax}(\chi)=\hat a(\chi)\hat x(\chi)\ne 0$ and hence $\hat a(\chi)\ne 0$ and
$\hat x(\chi)\ne 0$. Since $Supp(x)\subset null(a)$, $\hat a(\chi)\ne 0$ 
implies $\hat b(\chi)\ne 0$ for some $b\in ann(x)$. However,
$\hat{bx}(\chi)=\hat b(\chi)\hat x(\chi)\ne 0$, a contradiction.
\end{example}

Let $A'$ denote the Banach space dual to $A$. 
Setting $a\cdot x(\cdot)=x(a\cdot)$ for any $a\in A$, $x\in A'$, we have that
$A'$ is an $A$-module. This example is especially important because 
of the following result that connects the notions of spectral synthesis and synthesis
with respect to  $A$-modules.
 
\begin{theorem}\label{ind}
For $a\in A$ the following conditions are equivalent.

$1)$ $a$ admits spectral  synthesis;

$2)$ $a$ admits synthesis with respect to  $A'$;

$3)$ $a$ admits synthesis with respect to  any $A$-module.  
\end{theorem}
\begin{proof} $1)\Rightarrow 3)$. Let $M$ be an $A$-module and $x\in M$,
$Supp (x)\subset null (a)$. Let $ J = J(Supp (x))$,
then since 
$null(ann(x))=Supp(x)\subset null(a)$ and $a$ admits spectral synthesis in 
$A$, we have $a\in J_a\subset  J\subset ann(x)$ and hence $a\cdot x=0$.

$3)\Rightarrow 2)$ is obvious.

$2)\Rightarrow 1)$. Assume that $a\notin J_a$. Then there is $x\in A'$
such that $x(J_a)=0$ and $x(a)\ne 0$. Since $a$ admits synthesis with 
respect to
$A'$, $Supp(x)$ is not a subset of $null(a)$. On the other hand, 
$J_a\subset ann(x)$ 
and $Supp(x)=null(ann(x))\subset null(J_a)=null(a)$, a contradiction.
\end{proof}

\begin{lemma}\label{l1}
Let $M$ be a Banach $A$-module. 
For $x\in M$,  $a\cdot x=0$ for any $a\in J(E)$ if and only if 
$Supp(x)\subset E$.
\end{lemma}  
\begin{proof} Assume that $Supp(x)\subset E$ and  
that $\hat a=0$ on a nbhd of $E$. Then there exists an open
set $F$ such that $null(a)\supset F\supset Supp(x)$. Clearly,
$J=ann(x)$ is an ideal in $A$ with $null(J)=Supp(x)$. Since
$F^c\cap Supp(x)=\emptyset$,
there exists $b\in J$ such that $\hat b=1$ on $F^c$ and therefore 
$\hat{ab}=\hat a$ implying $a\cdot x=ab\cdot x=0$. 

If $a\cdot x=0$ for any $a\in J(E)$ then $J(E)\subset ann(x)$ and
$E=null (J(E))\supset null(ann(x))=Supp(x)$.
\end{proof}

As usually by $L_a$ we denote the operators of the ``left'' multiplication 
by $a \in A$,
acting in an $A$-module.
\begin{lemma}\label{l2}
Let $M$ be a Banach $A$-module and let $N$ be the closed submodule
generated by $x\in M$. Then $\sigma(L_a|_N)=\hat a(Supp(x))$.
\end{lemma}
\begin{proof}
Let $\tilde J=J(Supp(x))$.
Since, by Lemma~\ref{l1}, $b\cdot x=0$ for any $b\in \tilde J$, $N$ is also
an $A/\tilde J$-module with the action 
$(c+\tilde J)\cdot y=c\cdot y$. $X_{A/\tilde J}$ can be identified
with $\{\chi\in X_A:\chi(\tilde J)=0\}=Supp(x)$ and therefore
$\sigma(a+\tilde J)=\hat a(Supp(x))$. If $\lambda\in\hat a(Supp(x))^c$ then
$a-\lambda e+\tilde J$ is invertible in $A/\tilde J$ and so is the operator 
$(L_a-\lambda I)|_N$. This shows that 
$\sigma(L_a|_N)\subset \hat a(Supp(x))$.

Assume now that $(L_a-\hat a(\chi)I)|_N$, for some $\chi\in Supp(x)$, is an invertible
operator. Then there exists $c\in A$ such that 
$(a-\hat a(\chi))cb\cdot x=b\cdot x$ for any $b\in A$ and therefore 
$((a-\hat a(\chi))c-1)A\subset ann(x)$. Since 
$\chi\in Supp(x)=null(ann(x))$, we have $\chi(A)=0$. A contradiction. 
\end{proof}
\begin{theorem}\label{limit}
Let $M$ be a Banach $A$-module. For $x\in M$ and $a\in A$, $Supp(x)\subset
null(a)$ if and only if  $||a^nx||^{1/n}\to 0$ as $n\to\infty$.
\end{theorem}

\begin{proof}
Consider the Banach algebra $V=A/J_a$. The element $a+J_a$ is quasi-nilpotent
in $V$. In fact, if $\chi$ is a character of $A/J_a$ then $\rho(\chi)$, defined
by $\rho(\chi)(b):=\chi(b+J_a)$, $b\in A$, is a character of $A$. Therefore,
there exists $\tau\in X_A$ such that 
$\rho(\chi)(b)=\tau(b)$. This gives $\tau(b)=0$ for any $b\in J_a$ and hence
$\tau\in null(J_a)=null(a)$ and $\chi(a+J_a)=\tau(a)=0$.   

By Lemma~\ref{l1}, $J_a\cdot x=0$ whence $J_a\cdot N=0$, $N$ being
the closed submodule generated by $x$. $N$ becomes a Banach $V$-module 
by setting
$(b+J_a)\cdot y:=b\cdot y$, for $y\in N$ and $b\in A$.
Moreover,
$$||a^n\cdot x||^{1/n}=||(a+J_a)^n\cdot x||^{1/n}\leq ||(a+J_a)^{n}||^{1/n}
||x||^{1/n}\to 0, n\to\infty.$$

The converse follows immediately from Lemma~\ref{l2}.
\end{proof}

Let us also mention the module version of the "global" synthesis.
Let $M$ be a Banach $A$-module and $E$ be a closed subset of $X_A$.
$E$ is called a {\it set of synthesis over} $M$  ({\it synthetic over} $M$) 
if $a\cdot x=0$ for any 
$x\in M$ and $a\in A$ such that $Supp(x)\subset E\subset null(a)$.
Clearly, if, for $a\in A$, $null(a)$ is synthetic over $M$ then $a$ admits
synthesis with respect to  $M$.

\section{Modules over tensor algebras and lin\-ear operator equ\-ations}\label{tensor}
Let $X$ and $Y$ be compact Hausdorff spaces and consider
the projective tensor product $V(X,Y)=C(X)\hat\otimes C(Y)$.
Recall that $V(X,Y)$ (the Varopoulos algebra) 
consists of all functions 
$F\in C(X\times Y)$ which admit a representation 
\begin{equation}\label{equa}
F(x,y)=\sum_{i=1}^{\infty}f_i(x)g_i(y),
\end{equation}
 where $f_i\in C(X)$, $g_i\in
C(Y)$ and $$\sum_{i=1}^{\infty}||f_i||_{C(X)}||g_i||_{C(Y)}<\infty.$$ 
$V(X,Y)$ is a Banach algebra  with the norm
$$||F||_V=\inf\sum_{i=1}^{\infty}||f_i||_{C(X)}||g_i||_{C(Y)},$$
where $\inf$ is taken over all  representations of $F$ in the above form 
(see \cite{var}).
We note that $V(X,Y)$ is a semi-simple  regular Banach algebra with spectrum
$X\times Y$.

Any element of $V(X,Y)'$ can be identified with a 
bounded bilinear form $B(f,g)=\langle B,f\otimes g\rangle$ on
$C(X)\times C(Y)$ (a bimeasure, in short).

Let $M(X)$ denote the space of finite Borel measures on $X$.
 For $\mu\in M(X)$, $\nu\in M(Y)$, set $H_1=L_2(X,\mu)$, $H_2=L_2(Y,\nu)$.
Then $V(X,Y)\subset V^{\infty}(X,Y)\subset \Gamma(X,Y)$. Note that,
for ${\mathcal F}\subset V(X,Y)$, $\text{null }{\mathcal F}$ coincides
with $\displaystyle \cap_{F\in{\mathcal F}}F^{-1}(0)$.

By $M_f$, $M_g$ we denote the multiplication operators in $H_1$, $H_2$ by functions $f(x)$, $g(y)$ respectively.
 
Setting, for $F(x,y)=\sum_{n=1}^{\infty}f_n(x) g_n(y)\in V(X,Y)$ and 
$T\in B(H_1,H_2)$, 
\begin{equation}\label{act}
F\cdot T=\sum_{n=1}^{\infty}M_{g_n}TM_{f_n}
\end{equation}
we obtain a $V(X,Y)$-module structure on $B(H_1,H_2)$.
So, for $T\in B(H_1,H_2)$, we have  $Supp(T)=\cap\ null (F)$, the intersection 
being taken over all functions $F\in V(X,Y)$ such that $F\cdot T=0$. 
Now we compare $Supp(T)$ with the
``inner'' definitions of a support introduced in Section~\ref{pseudo}.

\begin{lemma}\label{ll2}
For $F\in V(X,Y)$ and $T\in B(H_1,H_2)$, if $F\cdot T=0$ then
$T$ is supported in $null(F)$.
\end{lemma}
\begin{proof}
Take Borel sets
$U\subset X$, $V\subset Y$ 
such that $(U\times V)\cap null(F)=\emptyset$ 
and consider $Q_VTP_U\in B(L_2(U,\mu), L_2(V,\nu))$. Then 
$\chi_U(x)\chi_V(y)F(x,y)\ne 0$ on $U\times V$ and if
$\Psi$ denote the set $\{FG\mid G\in \Gamma(U,V)\}\subset\Gamma(U,V)$,
 we have $\text{null }\Psi\cong\emptyset$. By \cite[Corollary~4.3]{sht}, $\Psi$ is 
dense in $\Gamma(U,V)$. As
$$0=\langle F\cdot Q_VTP_U,G\rangle=\langle Q_VTP_U,FG\rangle, \quad 
G\in\Gamma(U,V),$$
we obtain $Q_VTP_U=0$ and  therefore $null(F)$ supports $T$.
\end{proof}
\begin{proposition}\label{small}
$Supp(T)$ is the smallest closed set which supports the operator $T$.
\end{proposition}
\begin{proof} Set $E=Supp(T)$.
We show first that $E$ supports $T$. By Lemma~\ref{l1},
$F\cdot T=0$ for any $F\in V(X,Y)$ vanishing on a nbhd of $E$. 
Therefore, by Lemma~\ref{ll2}, $null(F)$ supports $T$ for each
$F\in J(E)$. 
Since $E=null (J(E))$, $E$ supports $T$.

Let $W\subset X\times Y$ be a closed set supporting $T$. By \cite[Theorem~4.3]{sht},
 given $F\in J(W)$, $\langle T,FG\rangle=0$ for any
$G\in \Gamma(X,Y)$ and hence 
$F\cdot T=0$. Applying now Lemma~\ref{l1}, we obtain $Supp(T)\subset W$, showing that 
$Supp(T)$ is the smallest closed set supporting $T$. 
\end{proof}

By the proposition we have therefore $\supp T \subset^M Supp(T)=\text{supp}_A T$.

Let $K_F$ denote the set of all operators $T$  satisfying   the condition
$F\cdot T=0$, that is $K_F = ann (F)$ in $B(H_1,H_2)$.
It coincides with the space of solutions of the linear operator equation
\begin{equation}\label{le}
\sum_{n=1}^{\infty}M_{g_n}TM_{f_n}=0,
\end{equation}
where $F(x,y)=\sum_{n=1}^{\infty}f_n(x)g_n(y)$.

It follows from Lemma~\ref{ll2}  and Proposition \ref{small} 
 that $Supp(K_F)\subset null(F)$.
\begin{proposition}\label{ref}
$$\text{\rm supp } K_F=Supp(K_F)=null(F).$$ 
\end{proposition}
\begin{proof}
It suffices to show that $null(F)\subset \supp K_F$.
Let $P=P(null(F))$ be the set of all pseudo-integral operators 
$T_{\sigma}$ with $supp(\sigma)\subset null(F)$ (see \cite[Section~1.5]{arv}). 
It follows from
\cite[Section~2.2]{arv} that $\supp P=null(F)$. On the other hand it is easy to see that
$$(F\cdot T_{\sigma}u,v)=\int\!\!\!\int F(x,y)u(x)v(y)d\sigma(x,y)=0,$$
for any $u\in H_1$, $v\in H_2$. Hence, $P\subset K_F$ and 
$null(F)=\supp P\subset \supp K_F$. 
\end{proof}
\begin{remark}\rm
The result can be proved without the use of pseudo-integral operators (see the proof
of a more general result, Proposition~\ref{min}, below).
\end{remark}
We say that $F\in V(X,Y)$ is {\it operator synthetic} with
respect to $(\mu,\nu)$ (we also write $(\mu,\nu)$-synthetic or
operator synthetic if $\mu$, $\nu$ are fixed) if it is synthetic
with respect to  the $V(X,Y)$-module $B(H_1,H_2)$.

The following proposition can be considered as a local version of Theorem~6.1
from [ShT]. Unlike the latter
it is "two-sided". 

\begin{proposition}\label{indiv}
$F\in V(X,Y)$ admits spectral synthesis if and only if  it is operator
synthetic for any choice of finite measures on $X$, $Y$.
\end{proposition}
\begin{proof}
The necessity follows from Theorem~\ref{ind}. To prove the sufficiency it is
enough to show that $F$ admits synthesis in $V(X,Y)'$.
Assume that synthesis fails for $F$. Then we  can find a bimeasure
$B\in V(X,Y)'$ such that $Supp (B)\subset null (F)$ and 
$F\cdot B\ne 0$. 
By the Grothendieck theorem \cite{grothendieck} there exist measures 
$\mu\in M(X)$, $\nu\in M(Y)$ 
such that    $$|\langle B,f\otimes g\rangle|\leq C||f||_{L_2(X,\mu)}
||g||_{L_2(Y,\nu)}$$ for any $f\in C(X)$ and $g\in C(Y)$. Thus
there exists an operator $T\in B(L_2(X,\mu),L_2(Y,\nu))$ such that
$$\langle T,\Psi\rangle=\langle B,\Psi\rangle,\quad \Psi\in V(X,Y),$$
where in the left hand side we used the inclusion $V(X,Y)\subset\Gamma(X,Y)$.
 Therefore, $F\cdot T\ne 0$.
We will get a contradiction if  prove that $Supp(T)\subset null(F)$.
For $G\in J_F$, $v\in C(X)$ and $w\in C(Y)$, we have
$$(G\cdot Tv,w)=\langle G\cdot B, v\otimes w\rangle=0,$$ 
implying $G\cdot T=0$ and $Supp(T)\subset null(G)$.
Since $null (J_F)=null(F)$, $Supp(T)\subset null(F)$.
\end{proof}
\begin{theorem}\label{spec}
$F\in V(X,Y)$ is operator synthetic if and only if 
$K_F={\mathfrak M}_{max}(null(F))$.
\end{theorem}
\begin{proof}
Assume $F\in V(X,Y)$ is operator synthetic with respect to $(\mu,\nu)$. 
We have $F\cdot T=0$ for each $T$ such that $Supp(T)\subset null(F)$.
As $Supp(T)$ is the smallest closed set which supports $T$, we have
that each $T\in {\mathfrak M}_{max}(null(F))$ is a solution of the equation
$F\cdot T=0$, i.e., $T\in K_F$. Conversely, if $T\in K_F$ then, by 
 Lemma~\ref{ll2},
$T$ is supported in $null(F)$.

Assume now that ${\mathfrak M}_{max}(null(F))=K_F$, but $F$ is not synthetic 
with respect to $B(H_1,H_2)$.
Then there exists $T$ with  $Supp(T)\subset null(F)$ and therefore 
$T\in {\mathfrak M}_{max}(null(F))$, such that $F\cdot T\ne 0$. A contradiction.
\end{proof}

For $F\in V(X,Y)$, let $\Delta_F$ denote the multiplication 
operator $X\mapsto F\cdot X$ on $B(H_1,H_2)$.  Then $K_F=\ker\Delta_F$. Let
${\mathcal E}_{\Delta_F}(0)$ be the  $0$-spectral  subspace of 
$\Delta_F$:
$${\mathcal E}_{\Delta_F}(0)=\{T\in B(H_1,H_2)\mid ||\Delta_F^n(T)||^{1/n}\to 0, n\to\infty\}.$$ 

Recall that if ${\mathcal L}$ is a subspace in $B(H_1,H_2)$
then its {\it reflexive hull}, $\text{Ref }{\mathcal L}$, is the set of all 
operators $A$
such that $Ax\in\overline{{\mathcal L}x}$ for any $x\in H_1$. ${\mathcal L}$ is
said to be {\it reflexive} if $\text{Ref }{\mathcal L}={\mathcal L}$.
\begin{proposition}\label{ref1}
$$\text{Ref } K_F={\mathfrak M}_{max}(null(F))={\mathcal E}_{\Delta_F}(0).$$
\end{proposition}
\begin{proof}
It is a standard fact (see \cite{arv} or \cite{sht}) that for arbitrary 
$L^{\infty}(X,\mu)\times L^{\infty}(Y,\nu)$-bimodule $G$ the space $\text{Ref }G$
consists of operators supported by $\text{supp }(G)$. Hence the first equality follows
from Proposition~\ref{ref}. By Theorem~\ref{limit}, ${\mathcal E}_{\Delta_F}(0)$ consists of all operators $T$ such that
$Supp(T)\subset null(F)$. But, by  Proposition~\ref{small}, this 
condition
is  equivalent  to $\supp T\subset null(F)$. Hence 
${\mathcal E}_{\Delta_F}(0)={\mathfrak M}_{max}(null(F))$.
\end{proof}
\begin{cor}\label{co}
The following are equivalent:

(a) $F\in V(X,Y)$ is operator synthetic;

(b) the solution space of the equation $F\cdot T=0$ is reflexive;

(c) $\ker\Delta_F={\mathcal E}_{\Delta_F}(0)$.
\end{cor}

\begin{proof}
(a)$\Rightarrow$ (b). It is easily seen that ${\mathfrak M}_{max}(null(F))$ 
is reflexive.
The implication follows now from Theorem~\ref{spec}.

 (b)$\Rightarrow$ (c). Follows from the equality
$\text{Ref }\ker\Delta_F={\mathcal E}_{\Delta_F}(0)$ which is due to Proposition~\ref{ref1}.

(c)$\Rightarrow$(a).
Let $T\in B(H_1,H_2)$ be an operator supported in $null(F)$. 
By Proposition~\ref{small},
$Supp(T)\subset null(F)$ and, by Theorem~\ref{limit}, 
$T$ is in ${\mathcal E}_{\Delta_F}(0)$ and therefore in $\ker\Delta_F$.
Hence, ${\mathfrak M}_{max}(null(F))\subset\ker\Delta_F$. The reverse
inclusion follows from 
Lemma \ref{ll2}. Thus  ${\mathfrak M}_{max}(null(F))=\ker\Delta_F$.
 The statement now follows from 
Theorem~\ref{spec}.
\end{proof}
Proposition~\ref{indiv} and Corollary~\ref{co} reduce the problem of 
verification of individual synthesis to a purely operator problem. 
The comparison of 0-spectral subspace and kernels for multiplication operators
will be one of the main topics in the further sections.

\section{Equations of more general form. Relations to ``glo\-bal'' operator synthesis}

The study of the  equations of the form (\ref{le}) is a part of the general 
theory 
of linear operator equations
\begin{equation}\label{genform}
\Delta(X)=\sum_{k\in K} B_kXA_k=0
\end{equation}
where $\{A_k\}_{k\in K}={\mathbb A}$, $\{B_k\}_{k\in K}={\mathbb B}$ are finite or countable
families of operators.

If ${\mathbb A}$, ${\mathbb B}$ are commutative families of normal
operators one says that (\ref{genform}) is a linear operator equation with 
normal coefficients. 
Among them equations (\ref{le}) correspond to those  whose coefficients satisfy 
the restriction
\begin{equation}\label{strongcoef}
\sum_{k\in K}||A_k||^2<\infty,\quad\sum_{k\in K}||B_k||^2<\infty.
\end{equation}
Indeed, realizing all $A_k$ and  $B_k$ as multiplication operators
by continuous functions $f_k$, $g_k$ on $L_2(X,\mu)$, $L_2(Y,\nu)$,
(\ref{genform}) can be rewritten in a form (\ref{le}); clearly 
$F(x,y)=\sum_{k\in K} f_k(x)g_k(y)\in V(X,Y)$.

It is more convenient sometimes to choose ``spectral'' realization of 
coefficient families.
Let $\sigma({\mathbb A})$, $\sigma({\mathbb B})$ be the maximal ideal spaces
of the unital $C^*$-algebras generated by the families ${\mathbb A}$ and 
${\mathbb B}$ respectively. To any $t\in\sigma({\mathbb A})$
we associate  a sequence $\lambda(t)=(t(A_1), t(A_2),\ldots)\in l_2$;
the map $t\mapsto\lambda(t)$ is continuous and identifies $\sigma({\mathbb A})$
with a compact subset of $l_2$. Thus $C^*({\mathbb A})$ can be considered as 
$C(\sigma({\mathbb A}))$ and the operators $A_i$ correspond to
the coordinate functions on $l_2$ (restricted to $\sigma({\mathbb A})$).
In a similar way we realize $C^*({\mathbb B})$. The space
$B(H_1,H_2)$ becomes a $V(\sigma({\mathbb A}),\sigma({\mathbb B}))$-module
with respect to the operation 
$$(f\otimes g)\cdot T=f({\mathbb B})Tg({\mathbb A}).$$
In particular, $\Delta(T)=F\cdot T$, where 
$F(\lambda,\mu)=\sum_{k\in K}\lambda_k\mu_k$.

Let $E_{\mathbb A}(\cdot)$,  $E_{\mathbb B}(\cdot)$ be the spectral measures of
${\mathbb A}$ and ${\mathbb B}$ (on $\sigma({\mathbb A})$, 
$\sigma({\mathbb B})$). We say that an operator $T$ is 
{\it supported in} 
$U\subset\sigma({\mathbb A})\times\sigma({\mathbb B})$ if 
$$E_{\mathbb B}(\beta)TE_{\mathbb A}(\alpha)=0$$ for any
Borel sets $\alpha\subset\sigma({\mathbb A})$, 
$\beta\subset\sigma({\mathbb B})$ such that $(\alpha\times\beta)\cap U
=\emptyset$.
The following result directly follows from Proposition~\ref{ref1}
if the families ${\mathbb A}$, ${\mathbb B}$  have  cyclic vectors
(one only needs to realize $H_1$, $H_2$ as the $L_2$-spaces for scalar 
spectral measures of  ${\mathbb A}$, ${\mathbb B}$). In the general case
the proof is similar to the  proof of Proposition~\ref{ref1}.
\begin{proposition}\label{S}
Let $S=\{(\lambda,\mu)\in\sigma({\mathbb A})\times\sigma({\mathbb B})\mid
\sum_{k\in K}\lambda_k\mu_k=0\}$. Then ${\mathcal E}_{\Delta}(0)$ consists of all operators 
supported in $S$.
\end{proposition}
In our further study of spectral behavior of multiplication operators $\Delta$ the following estimate will be useful.
\begin{lemma}\label{diam}
If $X\in{\mathcal E}_{\Delta}(0)$ then
$$||\Delta(X)||\leq 
2(\sum_{k\in K}||B_k||^2)^{1/2}||X||\text{\rm diam }\sigma({\mathbb A}).$$
\end{lemma}
\begin{proof}
 Let $S=\{(\lambda,\mu)\in\sigma({\mathbb A})\times\sigma({\mathbb B})\mid
\sum_{k\in K}\lambda_k\mu_k=0\}$. 
By Proposition~\ref{S}, if $X\in{\mathcal E}_{\Delta}(0)$ then
$X$ is supported in $S$ and therefore $X=E_{\mathbb B}(\hat\beta)X$
%$$\Delta(X)=E_{\mathbb B}(\hat\beta)\Delta(X),$$
where $\displaystyle\hat\beta=(\cup
\beta)^c$ with the union taken over all 
relatively open  $\beta\subset\sigma({\mathbb B})$ such that 
$(\sigma({\mathbb A})\times\beta)\cap S=\emptyset$.
We can certainly assume $\hat\beta=\sigma({\mathbb B})$, for if not
we set  $B_k'=B_k|_{E_{\mathbb B}(\hat\beta)H_2}$ and replace
$\Delta$ by the operator $\Delta'$ with coefficients $\{A_k\}$, $\{B_k'\}$ 
(acting on $B(H_1, E_{\mathbb B}(\hat\beta)H_2)$). 

Let $\{\lambda_k\}_{k\in K}\in\sigma({\mathbb A})$. Then
\begin{eqnarray*}
\Delta(X)
=\sum_{k\in K}B_kX\lambda_k+
\sum_{k\in K}B_kX(A_k-\lambda_kI).
\end{eqnarray*}
For each $f\in H$, we have
\begin{eqnarray*}
||(\sum_{k\in K}B_kX(A_k-\lambda_kI)f||^2\leq
(\sum_{k\in K}||B_k||\cdot ||X||\cdot||(A_k-\lambda_kI)f||)^2\leq\\
\leq||X||^2\sum_{k\in K}||B_k||^2\sum_{k\in K}||(A_k-\lambda_kI)f||^2\leq\\
\leq ||X||^2\sum_{k\in K}||B_k||^2
(\sum_{k\in K}(A_k-\lambda_k I)^*(A_k-\lambda_kI)f,f)\leq \\
\leq||X||^2
\sum_{k\in K}||B_k||^2(\text{diam }\sigma({\mathbb A}))^2||f||^2,
\end{eqnarray*}
and therefore
$$||\sum_{k\in K}B_kX(A_k-\lambda_kI)||\leq
||X||(\sum_{k\in K}||B_k||^2)^{1/2}(\text{diam }\sigma({\mathbb A})).$$

To estimate the norm of the first summand, note  that
the spectrum $\displaystyle\sigma(\sum_{k\in K}B_k
\lambda_k)$ is equal
to $\{\sum_{k\in K}\mu_k\lambda_k\mid (\mu_k)_{k\in K}\in\sigma({\mathbb B})\}$.
By our assumption, for given 
$\mu=(\mu_k)_{k\in K}\in\sigma({\mathbb B})$,
there exists $\lambda(\mu)=(\lambda_k(\mu))_{k\in K}\in\sigma({\mathbb A}))$ 
such that
$(\lambda(\mu),\mu)\in S$. Therefore,
\begin{eqnarray*}
|\sum_{k\in K}\mu_k\lambda_k|=|\sum_{k\in K}\mu_k(\lambda_k-\lambda_k(\mu))|
\leq (\sum_{k\in K}|\mu_k|^2)^{1/2}
(\sum_{k\in K}|\lambda_k-\lambda_k(\mu)|^2)^{1/2}\leq\\
\leq (\sum_{k\in K}||B_k||^2)^{1/2}(\text{diam }\sigma({\mathbb A}).
\end{eqnarray*}
 Thus 
$$||\sum_{k\in K}B_k\lambda_k||\leq  
(\sum_{k\in K}||B_k||^2)^{1/2}(\text{diam }\sigma({\mathbb A})),$$
implying
$||\sum_{k\in K}B_kX\lambda_k||\leq
(\sum_{k\in K}||B_k||^2)^{1/2}(\text{diam }\sigma({\mathbb A}))||X||$, and
$$||\Delta(X)||\leq 2(\sum_{k\in K}||B_k||^2)^{1/2}||X||\text{diam }\sigma({\mathbb A}).$$
\end{proof}

The condition (\ref{strongcoef}) is not necessary for a linear operator 
equation (\ref{genform}) to have sense. In many situations one may only 
suppose that
\begin{equation}\label{coef}
\sum_{k\in K} B_kB_k^*<\infty, \quad \sum_{k\in K} A_k^*A_k<\infty,
\end{equation}
which means that the norms of partial sums are  bounded and  the series
strongly converge. For equations with normal coefficients one can realize
${\mathbb A}$ and ${\mathbb B}$  
as families of multiplication operators on
the spaces $H_1=L_2(X,\mu)$, $H_2=L_2(Y,\nu)$:  
$$A_ku(x)=f_k(x)u(x), \quad B_kv(y)=g_k(y)v(y).$$ By (\ref{coef}),
$$\sum_{k\in K} |f_k(x)|^2<\infty,\quad \sum_{k\in K}|g_k(y)|^2<\infty.$$
In other words we may write (\ref{genform}) as $F\cdot X=0$ with 
$F\in V^{\infty}(X,Y)$ if $V^{\infty}(X,Y)$-module structure in 
$B(H_1,H_2)$ is defined by (\ref{act}) (now the series converge strongly).
It is not difficult to see that $F(x,y)\Psi(x,y)\in\Gamma(X,Y)$, for each
$\Psi(x,y)\in\Gamma(X,Y)$, and 
$\langle F\cdot T,\Psi\rangle=\langle T,F\Psi\rangle$, showing  that
the action  is well-defined.  

Our next aim is to show that in terms of $V^{\infty}(X,Y)$-module structure one
can describe ${\mathfrak M}_{min}(E)$ for arbitrary pseudo-closed set
$E\subset X\times Y$.

\begin{proposition}\label{min}
$${\mathfrak M}_{min}(E)=\{T\in B(H_1,H_2)\mid F\cdot T=0, 
\text{ for any } F\in V^{\infty}(X,Y)
\text{ that vanishes on $E$}\}.$$
\end{proposition}
\begin{proof}
Let us denote the right hand side of the equality by ${\mathfrak M}_e(E)$. 
Recall that ${\mathfrak M}_{min}(E)$
can be  characterized by $E$-pairs of projections 
$(P,Q)$ (see Section~\ref{pseudo}). 
%$${\mathfrak M}_{min}(E)=\{ T\in B(H_1,H_2)\mid Q(T\otimes 1)P=0,
% (P,Q)\in \Lambda({\mathfrak M}_{min}(E))\}$$
Let $(P,Q)$ be such a pair  and $P=P(x)=(P_{ij}(x))$, $Q=Q(y)=(Q_{ij}(y))$,
$P_{ij}(x)\in L^{\infty}(X,\mu)$, $Q_{ij}(y)\in L^{\infty}(Y,\nu)$ be
the matrix representations with respect to a fixed basis in $l_2$. 
Since $P$, $Q$ are projections, 
$$\sum_{j}P_{ij}(x)P_{ij}^*(x)=\sum_{j}P_{ij}(x)P_{ji}(x)=P_{ii}(x)\in 
L^{\infty}(X,\mu)$$ and
$$\sum_{j}Q_{jk}^*(y)Q_{jk}(y)=\sum_{j}Q_{kj}(y)Q_{jk}(y)=Q_{kk}(y)\in 
L^{\infty}(Y,\nu)$$ so that
$F_{ik}(x,y)=\sum_j P_{ij}(x)Q_{jk}(y)\in  
V^{\infty}(X,Y)$. Moreover, each $F_{ik}$ vanishes on $E$. Therefore, 
assuming 
$T\in {\mathfrak M}_e(E)$ we obtain $F_{ik}\cdot T=0$, for any $i,k$,
implying $Q(1\otimes T)P=0$ and hence $T\in {\mathfrak M}_{min}(E)$. 

Conversely, if $T\in{\mathfrak M}_{min}(E)$ then $\langle T,\Psi\rangle=0$ for
any $\Psi\in\Phi(E)$. Therefore for 
any $F\in V^{\infty}(X,Y)$ such that 
$F$ vanishes on $E$  and any $\Psi\in\Gamma(X,Y)$ the 
following holds
$$\langle F\cdot T,\Psi\rangle=\langle T,F\Psi\rangle=0$$
which implies $F\cdot T=0$, i.e. $T\in{\mathfrak M}_e(E)$.
\end{proof}

\begin{cor}\label{cr}
If a pseudo-closed set $E\subset X\times Y$ is a set of operator
synthesis then any
operator supported in $E$ satisfies each operator equation 
$F\cdot T=0$ with $F\in V^{\infty}(X,Y)$ vanishing on $E$.
\end{cor}

%Next theorem is another characterization of sets of synthesis.
\begin{cor}\label{equation}
If a pseudo-closed set $E\subset X\times Y$ is a set of synthesis and
$F_1,F_2\in V^{\infty}(X,Y)$ such that $\text{\rm null }F_i\cong E$ 
then the  
corresponding
linear operator equations $F_1\cdot T=0$ and $F_2\cdot T=0$ are equivalent.
\end{cor}
\begin{proof}
Let $T$ be a solution of the
equation $F_1\cdot T=0$. In order to show 
that $F_2\cdot T=0$
  it is sufficient, by Proposition~\ref{min},  to 
show that 
$T\in{\mathfrak M}_{max}(E)$ ($=\mathfrak M_{min}(E)$ in our case).

Take $U\subset X$, 
$V\subset Y$ such that $(U\times V)\cap E\cong\emptyset$ and
consider the operator $Q_VTP_U$ in $B(L_2(U,\mu),L_2(V,\nu))$. We have 
$\chi_U(x)\chi_V(y)F_1(x,y)\ne 0$ m.a.e. on $U\times V$. Let $\Psi$ denote 
the set
$$\{F_1\cdot F\mid F\in\Gamma(U,V)\}\subset\Gamma(U,V).$$ 
Then $\text{null }\Psi\cong\emptyset$ and, by  \cite[Corollary~4.3]{sht},
$\Psi$ is dense in $\Gamma(U,V)$. As
$$0=\langle F_1\cdot Q_VTP_U,F\rangle=\langle Q_VTP_U,F_1\cdot F\rangle, \quad 
 F\in \Gamma(U,V),$$ we obtain $Q_VTP_U=0$.
\end{proof}

%{\bf Comments.} Ne poluchaetsya u menya dokazat' v obratnuyu storonu. Nado 
%umet' kak to priblizhat' $G$, $null G\supset E$, i $null G$ ne soderzhit 
%okrestnosti $E$  funktsijami $FG$, $null F=E$ ili chto-to v etom rode.
%Esli etogo utverzhdeniya u nas ne budet, to nado izmenit' kommentarii posle
%Corollary 3.3.

%VNIM A razve zdes' ne vse dokazano? KON

%Hotelos' dokazat' takoe, chto esli dlya funktsij $F_1$, $F_2$ takih chto 
%$null F_1\cong null F_2\cong E$ uravnenija $F_1\cdot T=0$, $F_2\cdot T=0$ 
%ekvivalentny, to $E$ mn-vo sinteza. Dlja etogo nado pokazat', chto esli
%$F=0$ na $E$, no $null F\ne E$ to $F\cdot T=0$.    

%Conversely, let  $F_1\cdot T=0$ $\Leftrightarrow$ $F_2\cdot T=0$ for any
%$F_1$, $F_2\in V^{\infty}(X,Y)$, $\text{null }F_i\cong E$ and take 
%$T\in{\mathfrak M}_{max}(E)$. Since $E$ is pseudo-closed, $E\cong
%(\cup_{i=1}^{\infty}\alpha_i\times\beta_i)^c$ for some measurable
%$\alpha_i\subset X$, $\beta_i\subset Y$. Clearly, 
%$E\cap(\alpha_i\times\beta_i)\cong\emptyset$ and therefore
%$P_{\beta_i}TQ_{\alpha_i}=0$. Consider 
%$\displaystyle F(x,y)=\sum_{i=1}^{\infty}\frac{1}{2^i}\chi_{\alpha_i}(x)
%\chi_{\beta_i}(y)$. It is easy to see that $F\in V^{\infty}(X,Y)$, 
%$\text{null }F\cong E$ and $F\cdot T=0$. Therefore $F\cdot T=0$ for any
%$F\in V^{\infty}(X,Y)$ such that $\text{null F}\cong E$ by the assumption.
%From this one easily gets that $F\cdot T=0$ for any
%$F\in V^{\infty}(X,Y)$ such that $F$ vanishes on $E$. 
%By Proposition~\ref{min}, $T\in {\mathfrak M}_{min}(E)$, finishing the proof.
\begin{remark}\rm\label{sys}
The result extends to systems of equations $F_1^i\cdot T=0$, $1\leq i\leq n$.
In this case we have that $\supp T\subset\text{null }F_1^i$, for any $i$, and therefore
$\supp T\subset E$ for 
 $E=\cap_{i=1}^n\text{null }F_1^i$.
\end{remark}

\begin{cor}
Let $f_i$, $g_i$, $1\leq i\leq n$, be Borel functions
on standard Borel spaces
$(X,\mu)$, $(Y,\nu)$.
If $T\in B(L_2(X,\mu),L_2(Y,\nu))$ satisfies  operator equations 
$$M_{f_i}T=TM_{g_i}, \quad  1\leq i\leq n,$$
then $F\cdot T=0$ for any 
$F\in V^{\infty}(X,Y)$ vanishing on 
$\{(x,y)\mid f_i(x)=g_i(y), 1\leq i\leq n\}$.
\end{cor}

\begin{proof} By \cite[Theorem~4.8]{sht}
the set $\{(x,y)\mid f_i(x)=g_i(y), 1\leq i\leq n\}$ is synthetic. 
As it was noticed in Remark~\ref{sys},
$T$ is supported in $\{(x,y)\mid f_i(x)=g_i(y), 
1\leq i\leq n\}$ 
and the statement now follows from Corollary~\ref{cr}.
\end{proof}

The result implies, in particular, the Fuglede-Putnam Theorem,
 a useful tool in the operator theory,
which states the equivalence of the relations
$AT=TB$ and $A^*T=TB^*$, where $A$, $B$ are normal bounded operators
on a  Hilbert space and $T$ is just a bounded one acting on the 
same space.

It is natural to ask if the Fuglede-Putnam Theorem
extends to the equations of the form
$\sum_{i=1}^nB_iTA_i = 0$ and $\sum_{i=1}^nB_i^*TA_i^* = 0$, where 
$\{A_i\}_{1\le i\le n}$ and $\{B_i\}_{1\le i\le n}$ are commutative 
families of normal operators. This question of Gary Weiss 
\cite{weiss1}
has been answered  negatively in \cite{sh}. The proof  in \cite{sh} (see also 
\cite{SST} where it
is written more transparently) was based 
exactly on the connection between the individual and global
operator synthesis and related to the Schwartz example of a non-synthetic set
in the operator version due to Arveson \cite{arv}). 

In what follows we will 
find various conditions providing the equivalence for the equations 
of this kind and for more general linear operator equations. One of the 
main tools will be reducing (under some assumptions)
the problem to equations in the space of Hilbert-Schmidt operators. So we 
begin with a general approach that relates the spaces of the solutions of linear 
equations in different topological vector spaces.

\section{Approximate inverse intertwinings}

Let ${\mathfrak X}$ and ${\mathfrak Y}$ be topological vector spaces,
$\Phi:{\mathfrak X}\to{\mathfrak Y}$ a continuous imbedding with dense range, 
and $S$ and $T$ operators acting in ${\mathfrak X}$ and ${\mathfrak Y}$, 
respectively, intertwined by the mapping $\Phi$: $T\Phi=\Phi S$.
We write in this case that we are given an intertwining triple 
(or just an intertwining) $(\Phi,S,T)$. 

A net of linear mappings $F_{\alpha}:{\mathfrak Y}\to{\mathfrak X}$ is
called an {\it approximate inverse intertwining} (AII) for the
intertwining $(\Phi,S,T)$ if
$F_{\alpha}\Phi\to 1_{\mathfrak X}$, $\Phi F_{\alpha}\to 1_{\mathfrak Y}$ and
$F_{\alpha}T-SF_{\alpha}\to 0_{\X}$ in the topology of simple convergence.
%Replacing the last condition  by  the boundedness of the net
%$\{F_{\alpha}T-SF_{\alpha}\}$ in the same topology we arrive at the definition 
%of an {\it approximate inverse semi-intertwining} (AISI).

Denote by $\Phi^{-1}$ the full inverse image under the mapping
$\Phi$: $\Phi^{-1}(M)=\{x\in {\mathfrak X}\mid \Phi(x)\in M\}$ for any 
$M\subset {\mathfrak Y}$ 
(non-necessarily $M\subset\Phi({\mathfrak X})$). As usually
the  image of a map $X$ is denoted by $\im X$.
\begin{theorem}\label{aiipr1}
If the intertwinings $(\Phi, S_i, T_i)$, $1\leq i\leq n$,  have a common
AII, then $$\Phi^{-1}(\sum_i \im T_i)\subset \overline{\sum_i \im S_i}.$$
\end{theorem}

\begin{proof}

If $\Phi x=\sum_i T_iy_i$, then 
\begin{eqnarray*}
&\displaystyle x=\lim_{\alpha} F_{\alpha}\Phi x=\lim_{\alpha}
F_{\alpha}\sum_iT_iy_i=\\
&\displaystyle =\lim_{\alpha}(\sum_i(F_{\alpha}T_i-S_iF_{\alpha})y_i+
\sum_iS_iF_{\alpha}y_i)=\lim_{\alpha}\sum_iS_iF_{\alpha}y_i\in \overline{\sum_i \im S_i}.
\end{eqnarray*}
\end{proof}
Let ${\mathcal H}$ be a Hilbert space 
equipped with the weak operator topology.
\begin{cor}\label{aiicor1}
If ${\mathfrak X}={\mathcal H}$ and $(\Phi,S,T)$ has an AII, then
$$\Phi(\ker S^*)\cap\im T=\{0\}.$$
\end{cor}
\begin{proof}
$$\Phi(\ker S^*)\cap\im T =  \Phi(\ker S^*\cap\Phi^{-1}(\im T))\subset
\Phi(\ker S^*\cap\overline{\im S})=\{0\}.$$
\end{proof}

In applications $\X$, $\Y$ will be Banach spaces (of operators) supplied with 
weak or weak$^*$ topology. Nevertheless AII's can be used to obtain some 
norm inequalities. Such a possibility is provided by the following result:

\begin{proposition}\label{ad}
Suppose that $\X$ is a Banach space with the weak topology.
If $G_{\lambda}:\Y\to\X$, $\lambda\in\Lambda$, is an AII for an intertwining 
$(\Phi,S,T)$ then there is another  AII, 
$\{F_{\alpha}\}_{\alpha\in{\mathfrak A}}$, satisfying the conditions
\begin{equation}\label{norm1}
||F_{\alpha}\Phi x-x||\to 0,\text{ for any }x\in\X,
\end{equation} 
and 
\begin{equation}\label{norm2}
||(F_{\alpha}T-SF_{\alpha})y||\to 0,\text{ for any }y\in\Y.
\end{equation}
\end{proposition}
\begin{proof}
Let ${\mathfrak A}$ be the set of all triples $\alpha=(E,\lambda,\varepsilon)$,
where $E$ is a finite subset of $\X$, $\lambda\in\Lambda$, $\varepsilon >0$.
Setting $\alpha_1<\alpha_2$ if $E_1\subset E_2$, $\lambda_1<\lambda_2$,
$\varepsilon_1<\varepsilon_2$ we invert ${\mathfrak A}$ into a directed set.

Fix $\alpha=(E,\lambda,\varepsilon)\in{\mathfrak A}$. We claim that there is
a convex combination, $F=F_{\alpha}$, of operators $G_{\mu}$ with
$\mu>\lambda$ such that $||F\Phi x-x||<\varepsilon$ for any $x\in E$. 

Indeed, let $x_{\mu}=G_{\mu}\Phi x-x$. By our assumption $x_{\mu}\to 0$ (weakly in
$\X$), for any $x\in\X$. Let $N=\text{card }E$ and let $\X^N$ be  the direct sum of $N$ copies
of $\X$. Then the net $\displaystyle e_{\mu}=\oplus_{x\in E}x_{\mu}$ tends
to $0$ weakly  in $\X^N$, whence there is a convex combination 
$e=\sum_{i=1}^n c_ie_{\mu_i}$ with $\mu_i>\lambda$ such that
$||e||<\varepsilon$. Now setting $F=\sum_{i=1}^nc_iG_{\mu_i}$ we prove the claim.

It is clear that the net $\{F_{\alpha}\}_{\alpha\in{\mathfrak A}}$ is an AII for
$(\Phi, S,T)$ and that (\ref{norm1}) is satisfied. To obtain
(\ref{norm2}) one should repeat the trick (clearly the property 
(\ref{norm1}) will be preserved).
\end{proof}
\begin{theorem}\label{aiipr2}
Let $\Phi$ intertwine pairs $S_i$, $T_i$ ($i=1,2$). Suppose that
$\X$ is a Banach space equipped with a weak topology and 
$||S_2x||\leq ||S_1x||$ for any $x\in{\mathfrak X}$. If 
$(\Phi,S_1,T_1)$ has AII then
$$T_1^{-1}(\im\Phi)\subset T_2^{-1}(\im\Phi)$$
and 
\begin{equation}\label{*}
||\Phi^{-1}T_2y||\leq||\Phi^{-1}T_1y||
\end{equation}
for any $y\in T_1^{-1}(\im\Phi)$.
\end{theorem}
\begin{proof}
By Proposition~\ref{ad} we can assume that $(\ref{norm1})$ and $(\ref{norm2})$ 
hold.
Let $y\in T_1^{-1}(\im\Phi)$. Thus, $T_1y=\Phi x_1$ for some 
$x_1\in{\mathfrak X}$. Hence,
$$\displaystyle x_1=\lim_{\alpha}F_{\alpha}\Phi x_1=\lim_{\alpha}F_{\alpha}
T_1y=\lim_{\alpha} ((F_{\alpha}T_1-S_1F_{\alpha})y+S_1F_{\alpha}y)=\lim_{\alpha}S_1F_{\alpha}y.$$

Since $\{S_1F_{\alpha}y\}$ is a Cauchy net and 
$||S_2F_{\alpha}y-S_2F_{\beta}y||\leq 
||S_1F_{\alpha}y-S_1F_{\beta}y||$, we have that 
$\{S_2F_{\alpha}y\}$ is Cauchy. Let $x_2=\lim_{\alpha}S_2F_{\alpha}y$. Then
$||x_2||\leq ||x_1||$ and
$$\Phi x_2=\lim_{\alpha}\Phi S_2F_{\alpha}y=\lim_{\alpha} T_2\Phi F_{\alpha}y=T_2y,$$
the convergence being in the weak topology.
This imply $y\in T_2^{-1}(\im\Phi)$ and
$$||\Phi^{-1}T_2y||=||x_2||\leq ||x_1||=||\Phi^{-1}T_1y||.$$
\end{proof}
The following result has some similarity to Theorem~\ref{aiipr2} but
it does not use $AII's$.
\begin{proposition}\label{intertw}
Let $\Phi:\H\to\Y$ intertwine   a normal operator $S$   with $T_1$ and its 
adjoint $S^*$ with $T_2$. Suppose that 
$\displaystyle\ker S\cap\overline{\Phi^{-1}(T_2\Y)}=\{0\}$. Then (\ref{*}) holds for any
$y\in\Y$ such that $T_iy\in\Phi\H$, $i=1,2$.
\end{proposition}
\begin{proof}
Note first that $T_1$ and $T_2$ commute. Indeed,
$(T_1T_2-T_2T_1)\Phi y=\Phi(SS^*-S^*S)y=0$; since
$\Phi\H$ is dense in $\Y$ the claim follows.

Let $U$ be a partially isometric operator such that $SU=S^*$, 
$U\H=\overline{S\H}=\overline{S^*\H}$. Let $T_1y=\Phi x_1$, $T_2y=\Phi x_2$.
We have to prove that $||x_2||\leq||x_1||$. For $x=x_2-Ux_1$, one has
$$\Phi Sx=\Phi Sx_2-\Phi S^*x_1=T_1\Phi x_2-T_2\Phi x_1= T_1T_2y-T_2T_1y=0,$$
and hence  $Sx=0$, $x\in\ker S$.

On the other hand, $\Phi x_2=T_2y\in T_2\Y$, $x_2\in\Phi^{-1}(T_2\Y)$,
$Ux_1\in\overline{S\H}\subset\overline{\Phi^{-1}(T_2\Y)}$, so that
$x\in \overline{\Phi^{-1}(T_2\Y)}\cap\ker S=\{0\}$. Hence $x_2=Ux_1$,
$||x_2||\leq||x_1||$.
\end{proof}
We return to AII's. The following result is an immediate consequence
of Theorem~\ref{aiipr2}.
\begin{cor}\label{aiicor2}
Suppose that $S$ is a normal operator on $\H$ and  
$(\Phi,S,T_1)$, 
$(\Phi, S^*,T_2)$ have approximate inverse interwtinings (non-necessarily 
coinciding). Then 
$$||\Phi^{-1}T_2y||=||\Phi^{-1}T_1y||$$ for any 
$y\in T_1^{-1}(\im\Phi)=T_2^{-1}(\im\Phi)$
and,
in particular, $\ker T_1=\ker T_2$.
\end{cor}
In many cases the verification that a net $\{F_{\alpha}\}$ is an AII
can be considerably simplified by using the following result.
\begin{proposition}\label{add2}
Let $(\Phi,S,T)$ be an intertwining triple and $\{F_{\alpha}\}$ a net
of operators from $\Y$ to $\X$.

(i) If $\X$ is a dual Banach space with the weak-* topology and $\{F_{\alpha}\}$
satisfies the conditions

\hspace{0.4cm} (a) $\Phi F_{\alpha}y\to y$, for any $y\in\Y$;

\hspace{0.4cm} (b) $\{F_{\alpha}\Phi x\}$ is bounded for any $x\in\X$;

\hspace{0.4cm} (c) $\{(F_{\alpha}T-SF_{\alpha})y\}$ is bounded for any $y\in\Y$,

\noindent
then $\{F_{\alpha}\}$ is an AII.

\vspace{0.2cm}

(ii) If $\Y$ is a Banach space with the weak topology and 
$\{F_{\alpha}\}$ satisfies the conditions

\hspace{0.4cm} (d) $F_{\alpha}\Phi x\to x$, for any $x\in\X$;

\hspace{0.4cm} (e) $\sup_{\alpha}||\Phi F_{\alpha}||<\infty$;

\hspace{0.4cm} (f) for any neighbourhood $U$ of $0$ in $\X$ there is 
$\delta > 0$ such that $(F_{\alpha} T - S F_{\alpha})y \in U$ for all $\alpha$, when 
$||y|| < \delta$,

\noindent
then $\{F_{\alpha}\}$ is an AII.
\end{proposition}
\begin{proof}
$(i)$ Let $x\in\X$. We have to prove that $F_{\alpha}\Phi x\to x$. 
Since the net $\{F_{\alpha}\Phi x\}$ is bounded it is precompact (in the chosen
topology of $\X$) so it suffices to show that $x$ is its only limit point.
But if $x_1$ is a limit point of $\{F_{\alpha}\Phi x\}$ then $\Phi x_1$
is a limit point of $\{\Phi F_{\alpha}\Phi x\}$ which tends to $\Phi x$.
So $\Phi x_1=\Phi x$, $x_1=x$.

The proof of the condition $(F_{\alpha}T-SF_{\alpha})y\to 0$ is similar.

$(ii)$ The uniform boundedness permits us to prove the limit condition
$\Phi F_{\alpha}y\to y$ and $(F_{\alpha}T-SF_{\alpha})y\to 0$ on a dense subset.
But for $y\in\Phi\X$ they evidently follow from $(d)$.  
\end{proof}

In general a net $\{F_{\alpha}\}$ satisfying the conditions $(a)$, $(b)$, $(c)$ of
Proposi\-tion~\ref{add2} is called an 
{\it approximate inverse semi-intertwining} (AIS).

 Denote by ${\mathfrak X}^*$ the space of continuous antilinear functionals 
on ${\mathfrak X}$, endowed with the weak-* topology (in particular, 
${\mathcal H}^*={\mathcal H}$). The adjoint operators  (on $\X^*$ or between
$\X^*$ and $\Y^*$) are defined in the usual way. In particular, the
adjoint of an operator on ${\mathcal H}$ has the usual meaning.

It is not difficult to see that if $\{F_{\alpha}\}$ is an AII for 
$(\Phi,S,T)$ then $\{F_{\alpha}^*\}$ is an AII for $(\Phi^*, T^*, S^*)$.

Let ${\Phi}: {\mathcal H} \to {\mathfrak Y}$ intertwine operators 
$S, S^*$ with $T_1,T_2$. Let 
$\{F_{\alpha}\}:{\mathfrak Y}\to{\mathcal H}$ be an AII for the 
intertwining $(\Phi,S,T_1)$. It is called
 a {\it $*$-approximate inverse intertwining} ($*$-AII) for the ordered
 pair ($(\Phi,S,T_1)$, $(\Phi,S^*,T_2)$)
  if $\{F_{\alpha}^*F_{\alpha}\}$ is an AII for $(\Phi\Phi^*,T_1^*,T_2)$.

A {\it $*$-approximate semi-intertwining} ($*$-AIS) is defined in a similar way: it is an AIS $\{F_{\alpha}\}$
such that $\{F_{\alpha}^*F_{\alpha}\}$ is an AIS for
$(\Phi\Phi^*,T_1^*,T_2)$. Since $\H$ is reflexive, $\{F_{\alpha}\}$ is in fact 
an AII by Proposition~\ref{add2}.
\begin{cor}
If $\Y$ is a Banach space with the weak topology then any $*$-AIS is a $*$-AII. \end{cor}
\begin{proof}
Follows from Proposition~\ref{add2}$(i)$ (with $\Y^*$ as the space $\X$).
\end{proof}

\begin{theorem}\label{aiipr3}
(i) If the pair $((\Phi, S,T_1)$, $(\Phi,S^*,T_2))$ has $*$-AIS, then
$$(\im T_1)\cap T_2^{-1}(\Phi\Phi^*({\mathfrak Y}^*))\subset
\Phi({\mathcal H}).$$ 
(ii) If $((\Phi, S,T_1)$, $(\Phi,S^*,T_2))$ has $*$-AII  then
$$||\Phi^{-1}(T_1y)||^2=\langle (\Phi\Phi^*)^{-1}(T_2T_1y), y\rangle$$
for any $y\in (T_2T_1)^{-1}(\Phi\Phi^*(\Y^*))$.
\end{theorem}
\begin{proof}
Let  $\hat y\in(\im T_1)\cap T_2^{-1}(\Phi\Phi^*({\mathfrak Y}^*))$. Then
$\hat y=T_1y$ and $T_2\hat y=\Phi\Phi^*z$, for some $z\in\Y^*$ and $y\in\Y$.
Let $\{F_{\alpha}\}$ be a $*$-AIS for  the pair $((\Phi, S,T_1)$, $(\Phi,S^*,T_2))$.
Then $\{F_{\alpha}^*F_{\alpha}T_2\hat y-T_1^*F_{\alpha}^*
F_{\alpha}\hat y\}$ is bounded.
Because $F_{\alpha}^*F_{\alpha}T_2\hat y=F_{\alpha}^*F_{\alpha}
\Phi\Phi^*z\to z$, we obtain also boundedness of the  net $\{||F_{\alpha}\hat y||^2\}$:
$$||F_{\alpha}\hat y||^2=(F_{\alpha}^*F_{\alpha}\hat y,T_1y)=(T_1^*F_{\alpha}^*F_{\alpha}\hat y,y).$$
Thus there exists a subnet $\{F_{n(\alpha)}\hat y\}$ 
converging
weakly to a vector $h\in{\mathcal H}$. This gives us
$$\hat y=\lim_{\alpha}\Phi F_{n(\alpha)}\hat y=\Phi h\in\Phi({\mathcal H}),$$
and the first statement of the theorem.

Let now $\{F_{\alpha}\}$ be a $*$-AII.
Then the net $\{F_{\alpha}^*F_{\alpha}T_2T_1y-
T_1^*F_{\alpha}^*F_{\alpha}T_1y\}$ converges to zero in the 
weak$^*$-topology in
${\mathfrak Y}^*$. By the previous arguments we have that 
$F_{\alpha}^*F_{\alpha}T_2T_1y\to z$, 
where $T_2T_1y=\Phi\Phi^*z$, implying  
$\{T_1^*F_{\alpha}^*F_{\alpha}T_1y\}\to z$  and
$$\langle T_1^*F_{\alpha}^*F_{\alpha}T_1y,y\rangle\to \langle z,y\rangle.$$
On the other hand,
by the first statement, $T_1y=\Phi h$ for some $h\in \H$ and 
$$\langle T_1^*F_{\alpha}^*F_{\alpha}T_1y,y\rangle=\langle F_{\alpha}^*F_{\alpha}T_1y,T_1y\rangle=\langle
F_{\alpha}^*F_{\alpha}\Phi h,\Phi h\rangle=(\Phi^*F_{\alpha}^*F_{\alpha}\Phi h, h)\to ||h||^2.$$
One can see the convergence here in the following way. Since  
$\{F_{\alpha}\Phi h\}$, $h\in\H$, is bounded, the net
$\{\Phi^*F_{\alpha}^*F_{\alpha}\Phi h\}$ is bounded and therefore  precompact in the
weak topology of $\H$. Moreover, as 
$\Phi\Phi^*F_{\alpha}^*F_{\alpha}\Phi h\to \Phi h$ and $\Phi$ is injective, 
we have that the only limit point of $\{\Phi^*F_{\alpha}^*F_{\alpha}\Phi h\}$ 
is $h$.
 
Finally we obtain
$$||\Phi^{-1}(T_1y)||^2=||h||^2=\langle z,y\rangle=\langle (\Phi\Phi^*)^{-1}(T_2T_1y), y\rangle.$$  
\end{proof}

\begin{cor}\label{aiicor3}
If $((\Phi,S,T_1)$, $(\Phi,S^*,T_2))$ has $*$-AIS then $\im T_1\cap\ker T_2=\{0\}$.
\end{cor}
\begin{proof}
Clearly, $\ker T_2\subset T_2^{-1}(\Phi\Phi^*(\Y^*))$. Therefore,
by Theorem~\ref{aiipr3}, $\im T_1\cap\ker T_2\subset\Phi({\mathcal H})$ 
and, by Corollary~\ref{aiicor1},
$$\im T_1\cap\ker T_2=\im T_1\cap\ker T_2\cap\Phi({\mathcal H})=\im T_1\cap\Phi(\ker S^*)=0.$$
\end{proof}

\section{AII for inclusions of symmetrically normed ideals and multiplication
operators}
%The purpose of this section is to apply the technique of 
%approximate intertwinings to study of linear 
%operator equations
%We are interested in  linear
%operator equations of the form 
%$$\sum_{k\in K}B_kXA_k=0,$$
%where $\{A_k\}_{k\in K}$ and $\{B_k\}_{k\in K}$ are
%commutative families of normal operators acting on $H$
%such that 
%$$\sum_{k\in K}||A_k||\cdot ||B_k||<\infty.$$
%Before proceeding to this  we first have to  establish some properties
%of multiplication operators $X\mapsto\sum_{k\in K}B_kXA_k$ on
%$B(H)$ and its symmetrically normed ideals.
 
Here we apply the results of the previous chapter to multiplication 
operators on symmetrically normed ideals of operators on  Hilbert spaces. 

Let $H$  be a Hilbert space, $B(H)$ be the space of
bounded linear operators on $H$.  
For a symmetrically normed ideal $J$ we denote by $||\cdot||_{J}$ the 
associated norm. If $J=\S_p$, $1 \le p< \infty$, a Shatten-von Neumann ideal, 
we simply write $||\cdot ||_p$ instead of $||\cdot||_{\S_p}$.
 
Given $J$, we set
\begin{eqnarray*}
J^*=\{X\in B(H)\mid XY\in{\mathfrak S}_1, \text{ for any }Y\in J\},\\
J'=\{X\in B(H)\mid XY\in{\mathfrak S}_2, \text{ for any }Y\in J\},\\
\tilde J=\{X\in B(H)\mid XY\in J^*, \text{ for any }Y\in J\}.
\end{eqnarray*}
It is clear that in this way one obtains ideals of $B(H)$ which become
symmetrically normed if the norm of an operator $X$ in $J^*$ ($J'$ and 
$\tilde J$) 
is defined as the norm
of the mapping $Y\mapsto XY$ from $J$ to $\S_1$ (from $J$ to $\S_2$ and
from $J$ to $J^*$ respectively).
 
One can easily see that
\begin{eqnarray*}
B(H)^*= {\S}_1, \quad  {\S}_1^*=B(H), \quad {\S}_{\infty}^*=\S_1,\\
B(H)'=\S_2, \quad \S_1'=B(H),\quad
\S_{\infty}'=\S_2,\\
 \tilde B(H)=\S_1,\quad \tilde\S_1=B(H),\quad 
\tilde\S_{\infty}=\S_1,
\end{eqnarray*}
and if $J=\S_p$, $1<p\leq\infty$,  then
\begin{eqnarray*}
J^*={\mathfrak S}_{p/(p-1)},\quad J'=\left\{\begin{array}{ll}
{\mathfrak S}_{2p/(p-1)},& p>2,\\
B(H), &p\leq 2,\end{array}\right.
 \quad \tilde J=\left\{\begin{array}{ll}
{\mathfrak S}_{p/(p-2)},& p>2,\\
B(H), &p\leq 2\end{array}\right.
\end{eqnarray*}
(the equality for symmetrically normed  ideals assumes the equality of 
the norms).

Let $\{A_k\}_{k\in K}$ and $\{B_k\}_{k\in K}$ be  arbitrary families 
of operators (not necessarily
commuting or normal) acting 
on $H$ such that $$\sum_{k\in K}||A_k||\cdot ||B_k||<\infty.$$
Note that 
multiplying by constants $A_k\mapsto\lambda_kA_k$, 
$B_k\mapsto\lambda_k^{-1}B_k$ we may (and will) assume that
\begin{equation}\label{red}
\sum_{k\in K}||A_k||^2<\infty,\qquad \sum_{k\in K}||B_k||^{2}<\infty.
\end{equation}
It is easy to check that in this case the multiplication 
operator $\Delta:X\mapsto \sum_{k\in K}B_kXA_k$ is continuous on $B(H)$ and 
preserves
all symmetrically normed ideals.
We shall also denote by $\tilde\Delta$ the formal adjoint to $\Delta$:
$\tilde\Delta (X) = \sum_{k\in K}B_k^*XA_k^*$.
Note that $\tilde\Delta|_{\S_2}=(\Delta|_{\S_2})^*$ and
$\tilde\Delta|_{\S_1}=\Delta^*$.
For simplicity of notation, we write $\Delta_J$ and $\tilde\Delta_J$ 
instead of
$\Delta|_J$ and $\tilde\Delta|_J$ for a symmetrically normed ideal 
$J$ of $B(H)$.

If $J_1 \subset J_2$ then the natural inclusion $J_1 \hookrightarrow J_2$ 
will be denoted by 
$\Phi_{J_1,J_2}$. Clearly $\Phi_{J_1,J_2}$ intertwines $\Delta_{J_1}$ with 
$\Delta_{J_2}$. For brevity we will denote this intertwining triple by 
$(\Phi ,\Delta_{J_1}, \Delta_{J_2})$.

Now we should look for the approximate inverse intertwinings. They will
be constructed by means of increasing sequences $\{P_n\}$ 
of finite-dimensional projections and will have the form $F_n(X) = XP_n$. 
But for this the coefficient families of a multiplication operator must
satisfy some restrictions.

For any family $\{X_k\}_{k\in K}$ of operators and a finite-dimensional
projection $P$ we set
$$\varphi_P^{J}(\{X_k\}_{k\in K})=(\sum_{k\in K}||[X_k,P]||^2_J)^{1/2},$$
where $||\cdot||_J$ is the norm in the ideal $J$.

A family $\{X_k\}_{k\in K}$ is said to be {\it $J$-semidiagonal} if
there exists a sequence of projections $P_n$ of finite rank such that
$P_n\to^s 1$ and $\sup_n\varphi_{P_n}^{J}(\{X_k\}_{k\in K})<\infty$.
If $J=\S_p$, $1\leq p\leq\infty$, we write  simply {\it $p$-semidiagonal}.
It is clear that if $p_1 < p_2$ then each $p_1$-semidiagonal family is 
$p_2$-semidiagonal. In particular, 1-semidiagonality is the strongest of 
these conditions. Clearly, any finite family is $\S_{\infty}$-semidiagonal.

Examples of semidiagonal families will be discussed later on.

In the following theorem $J$ stands for either a separable symmetrically
normed ideal with the weak topology or for a symmetrically normed ideal, 
dual to a separable one, with the weak$^*$-topology.
In both cases $J^*$ is identified with the space of continuous antilinear 
functionals by means of the map
$f_Y(X)=tr(X^*Y)$. Under this correspondence 
 $\tilde\Delta_{J^*}=\Delta_J^*$.

%for  a given
%$s>1$, 
%$$\sum_{k\in K}||A_k||^s<\infty,\qquad \sum_{k\in K}||B_k||^{s/(s-1)}<\infty,$$
%($\sup_k||B_k||<\infty$ if $s=1$). In fact, take 
%$(||B_k||^{1/s}/||A_k||^{s/(s-1)})A_k$  instead of $A_k$ and  
%$(||A_k||^{s/(s-1)}/||B_k||^{1/s})B_k$  instead of $B_k$ in case $s>0$ and
%$||B_k||A_k$ and $(1/||B_k||)B_k$ instead of $A_k$ and $B_k$ if $s=1$.

\begin{theorem}\label{aiipr4}
Assume  that $\S_2\subset J$. 

(i) If $\displaystyle\{A_k\}_{k\in K}$ is $J'$-semidiagonal  then
there is an AII for $(\Phi,\Delta_{\S_2},\Delta_J)$.  

\vspace{0.1cm}

(ii) If $\displaystyle\{A_k\}_{k\in K}$ is $\tilde J$-semidiagonal  then
there is a $*$-AIS for 
$((\Phi,\Delta_{\S_2},\Delta_J)$, $(\Phi,\tilde\Delta_{\S_2},\tilde\Delta_J))$,
and, moreover, a $*$-AII if $J$ is separable with the weak topology.

\vspace{0.1cm}
Assume that $\S_1\subset J$.  

(iii) If $\displaystyle\{A_k\}_{k\in K}$ is $J^*$-semidiagonal  then
there exists an AIS for $(\Phi,\Delta_{\S_1},\Delta_J)$ in general and
an AII if $J$ is separable with the weak topology.
\end{theorem}

\begin{proof}
$(i)$ 
We define  $F_n:J\to \S_2$ by $F_n(X)=XP_n$, $X\in J$, where 
$\{P_n\}$ is a sequence of finite rank projections  such that 
$\sup_n\varphi_{P_n}^{J}(\{A_k\}_{k\in K})<\infty$ and $P_n\to^s 1$, 
$n\to\infty$.
Clearly, $F_n\Phi\to 1$,
and $\Phi F_n\to 1$. For $X\in J$ one can easily check the equality 
$$(\Delta_{\S_2}F_n-F_n\Delta_J)(X)=\sum_{k\in K}B_kX[A_k,P_n]$$
and
\begin{eqnarray*}
&&||(\Delta_{\S_2}F_n-F_n\Delta_J)(X)||_2\leq\sum_{k\in K}||B_kX||_J
||[A_k,P_n]||_{J'}\leq\\
&&\leq(\sum_{k\in K}||B_kX||_J^2)^{1/2}
\varphi_{P_n}^{J'}(\{A_k\}_{k\in K})\leq||X||_J
(\sum_{k\in K}||B_k||^2)^{1/2}\varphi_{P_n}^{J'}(\{A_k\}_{k\in K}),
\end{eqnarray*}
showing that  $\{F_n\}$ is an AIS. 
Since $\S_2$ is reflexive, $\{F_n\}$ is an AII  by Proposition~\ref{add2}$(i)$.
%We next claim that $(\Phi, \Delta_{\S_2},\Delta_J)$  has  an AII.
%In fact, since $\S_2$ is reflexive and $\{(\Delta_{\S_2}F_n-F_n\Delta_J)(X)\}$
%is bounded for each $X\in J$, we have that for every $X\in J$ there exists a 
%subsequence $\{(\Delta_{\S_2}F_{n_k}-F_{n_k}\Delta_J)(X)\}$ which converges 
%weakly. As
%$$\Phi(\Delta_{\S_2}F_{n_k}-F_{n_k}\Delta_J)(X))=\Delta_J\Phi F_{n_k}(X)-
%\Phi F_{n_k}\Delta_J(X)\to 0,$$
%we obtain $$(\Delta_{\S_2}F_{n_k}-F_{n_k}\Delta_J)(X)\to^w 0.$$
%It is clear now, that for each finite set $\Sigma\subset J$ one can find
%a subsequence $\{\Delta_{\S_2}F_{n_k}-F_{n_k}\Delta_J\}$
%converging weakly to $0$ on every $X\in\Sigma$.
%Therefore, given $\Sigma$ and a neighbourhood $\Omega$ of zero in the 
%weak topology there exists $n(\Omega,\Sigma)$ such that 
%$(\Delta_{\S_2}F_{n}-F_{n}\Delta_J)(X)\in \Omega$ for each $X\in \Sigma$.

%Let $\Lambda$ be the  set of all pairs $(\Sigma,\Omega)$, where
%$\Sigma$ is a finite subset of 
%$J$ and 
%$\Omega$ is a neighbourhood of zero in the weak topology, with ordering defined
%as follows:
%$(\Sigma_1,\Omega_1)>(\Sigma_2,\Omega_2)$ if
%$\Sigma_1\supset\Sigma_2$ and $\Omega_1\subset\Omega_2$.
%Consider now a net $\{F_{\lambda}\}_{\lambda\in\Lambda}$
%with $F_{\lambda}=F_{n(\Omega,\Sigma)}$ for $\lambda=(\Sigma,\Omega)$.
%Then  for each $X\in J$,
%$$\displaystyle\lim_{\lambda\in\Lambda}
%(\Delta_{\S_2}F_{\lambda}-F_{\lambda}\Delta_J)(X)=0.$$

$(ii)$
Define $F_n:J\to\S_2$ as before: $F_n(X)=XP_n$, $X\in J$,
where $\sup_n\varphi_{P_n}^{\tilde J}(\{X_k\}_{k\in K})<\infty$ and 
$P_n\to^s 1$, $n\to\infty$.

Similar arguments shows that 
\begin{eqnarray*} 
&&||\Delta_J^*F_n^*F_n(X)-
F_n^*F_n\tilde\Delta_J(X)||_{J^*}
\leq\sum_{k\in K}||B_k^*X[A_k^*,P_n]||_{J^*}\leq\\
&&\leq\sum_{k\in K}||B_k^*X||_J
||[A_k^*,P_n]||_{\tilde J}
\leq ||X||_J
(\sum_{k\in K}||B_k||^2)^{1/2}\varphi_{P_n}^{\tilde J}(\{A_k\}_{k\in K}).
\end{eqnarray*}
giving $||\Delta_J^*F_n^*F_n-
F_n^*F_n\tilde\Delta_J||<\infty$. Thus $\{F_n^*F_n\}$ is an AIS for  
$(\Phi\Phi^*,\Delta_J^*,\tilde\Delta_J)$ and if $J$ is supplied with the weak
topology it is even 
an AII by Proposition~\ref{add2}.

%If $J$ is separable we can choose a subsequence,
%$\{F_{n_k}^*F_{n_k}\}$, which will be an AII.

The sequence 
 $\{F_n\}$ is an  AIS for 
$(\Phi,\Delta_{\S_2},\Delta_J)$ and, by reflexivity of   $\S_2$,
 it is also an AII.
In fact,
\begin{eqnarray*}
&&||(\Delta_{\S_2}F_n-F_n\Delta_J)(X)||_2\leq\sum_{k\in K}||B_kX[A_k,P_n]||_2\leq\\
&&\leq C\sum_{k\in K}||B_kX[A_k,P_n]||_{J^*}\leq\sum_{k\in K}||B_kX||_J
||[A_k,P_n]||_{\tilde J}\leq\\
&&\leq (\sum_{k\in K}||B_kX||_J^2)^{1/2}
\varphi_{P_n}^{\tilde J}(\{A_k\}_{k\in K})\leq||X||_J
(\sum_{k\in K}||B_k||^2)^{1/2}\varphi_{P_n}^{\tilde J}(\{A_k\}_{k\in K}).
\end{eqnarray*} 
(we used the fact that $\S_2\subset J$ and therefore $J^*\subset\S_2$ so that
$||\cdot||_2\leq C||\cdot||_{J^*}$ for some constant $C$). 

%(esli $J\ne \S_p$, to v neravenstve mogut 
%byt' i konstanty ili eshche chto-nibud' bolee slozhnoe???)
%VNIM Vrode by, konstant ne budet: vklyuchenie idealov vlechet ner-vo dlya norm.
%No my na etu chast' posmotrim vnimatel'nee pri sleduyuschej redakcii,
%kogda u menya budet pered soboj kniga Gohberga-Krejna, ona
%v Londone est') KON

%Thus $\{F_n\}$ and therefore $\{F_{n_k}\}$ is an AIS for 
% $(\Phi,\Delta_{\S_2},\Delta_J)$. To find  a net which will be an AII for  
%$(\Phi,\Delta_{\S_2},\Delta_J)$ and a $*$-AIS  and a $*$-AII in case 
%$J\ne B(H)$ for 
%$(\Phi,\Delta_{\S_2},\Delta_J)$, $(\Phi,\tilde\Delta_{\S_2},\tilde\Delta_J)$ 
%we proceed as in (i). 

$(iii)$ In a similar way  one shows that
$F_n:J\to\S_1$, $F_n(X)=XP_n$, $X\in J$, is an AIS for the intertwining
$(\Phi,\Delta_{\S_1},\Delta_J)$. Moreover, 
$$||\Delta_{\S_1}F_n-F_n\Delta_J||\leq \sum_{k\in K}(||B_k||^2)^{1/2}
\varphi_{P_n}^{J^*}(\{A_k\}_{k\in K}).$$
Therefore, to prove that  $\{F_n\}$ is an AII for separable $J$ 
(endowed with the weak topology) it is sufficient
to prove that
$B_kX[A_k,P_n]\to^w 0$, as $n\to\infty$, for any $X\in J$.

Given $Z\in B(H)$, 
\begin{eqnarray*}
|\tr(ZB_kX[A_k,P_n])|=|\tr(ZB_kX(1-P_n)A_kP_n)-\tr(ZB_kXP_nA_k(1-P_n))|\leq\\
\leq |\tr(ZB_kX(1-P_n)A_kP_n)|+|\tr((1-P_n)ZB_kXP_nA_k(1-P_n))|\leq\\
\leq ||(ZB_kX)(1-P_n)||_J\cdot||(1-P_n)A_kP_n||_{J^*}+||(1-P_n)(ZB_kX)||_J
\cdot||(P_nA_k(1-P_n)||_{J^*}.
\end{eqnarray*}
 Since $\sup_n||(P_nA_k(1-P_n)||_{J^*}<\infty$, and 
$||ZB_kX(1-P_n)||_J$, $||(1-P_n)ZB_kX||_J\to 0$ as $n\to\infty$ 
if $J$ is separable, we have the statement.
\end{proof}

\begin{remark}\label{rem}\rm
As in the proof of $(i)$ ($(iii)$ respectively) one can  show that having
a finite number of multiplication operators $\Delta_i:B(H)\to B(H)$,
$\Delta_i(X)=\sum_{k\in K}B_k^iXA_k^i$ such that the family  
$\{A_k^i\}_{i=1,k\in K}^n$ is $J'$-semidiagonal ($\tilde J$-semidiagonal)
there exists a common AII for the intertwinings
$(\Phi, (\Delta_i)_{\S_2},(\Delta_i)_J)$, $1\leq i\leq n$ 
($(\Phi, (\Delta_i)_{\S_1},(\Delta_i)_J)$, $1\leq i\leq n$).
\end{remark}

\begin{cor}\label{aiicor4}
(i) If $\{A_k\}_{k\in K}$ is $1$-semidiagonal, then there exist
an AIS for $(\Phi,\Delta_{\S_1},\Delta)$, an AII 
for $(\Phi,\Delta_{\S_2},\Delta_{\infty})$, a $*$-AIS for 
$((\Phi,\Delta_{\S_2},\Delta)$,
$(\Phi,(\Delta_{\S_2})^*,\tilde\Delta))$ and a $*$-AII for 
the pair $((\Phi,\Delta_{\S_2},\Delta_{\infty})$, $(\Phi,(\Delta_{\S_2})^*,\tilde\Delta_{\infty}))$.

\vspace{0.1cm}

(ii) If $\{A_k\}_{k\in K}$ is $2$-semidiagonal, 
then there exists an AII for $(\Phi,\Delta_{\S_2},\Delta)$.

\vspace{0.1cm}

(iii) If $\{A_k\}_{k\in K}$  is $p/(p-1)$-semidiagonal, then 
there exists an AII for $(\Phi,\Delta_{\S_1},\Delta_{\S_p})$.

\vspace{0.1cm}

(iv) If $\{A_k\}_{k\in K}$ is $2p/(p-2)$-semidiagonal, then
there exists an AII for $(\Phi,\Delta_{\S_2},\Delta_{\S_p})$.

\vspace{0.1cm}

(v) If $\{A_k\}_{k\in K}$ is $p/(p-2)$-semidiagonal,
then there exists a $*$-AII for 
$((\Phi,\Delta_{\S_2},\Delta_{\S_p})$, 
$(\Phi,(\Delta_{\S_2})^*,\tilde\Delta_{\S_p}))$.
\end{cor}

Now we list some  examples of semidiagonal families.
\begin{proposition}
If in some basis the matrices of all the operators $A_k$ have all their
 nonzero 
entries on a finite number of diagonals (and $\sum_{k\in K}||A_k||^2<\infty$),
then the family $\{A_k\}_{k\in K}$ is $1$-semidiagonal.
\end{proposition}
\begin{proof}
Let $\{e_k\}$ be a basis satisfying the assumptions. Then $(A_ke_i,e_j)=0$ for
$|i-j|>n$, where $n$ is a positive integer. Let $P_m$ be the projection onto
the subspace generated by $e_1,\ldots, e_k$. One can easily see that for each 
$m$ and $k$ the rank of  the operator $[A_k,P_m]$ does not exceed $2n+1$ and
therefore
$$||[A_k,P_m]||_1\leq (2n+1)||[A_k,P_m]||\leq 2(2n+1)||A_k||,$$
$$\sup_m\varphi_{P_m}^{\S_1}(\{A_k\}_{k\in K})\leq 2(2n+1)
(\sum_{k\in K}||A_k||^2)^{1/2}<\infty.$$
\end{proof}
The simplest class of such   examples consists of finite families of weighted
shifts.

More generally one can consider operators with matrices whose
entries $a_{ij}$ sufficiently quickly decrease with $|i-j|\to\infty$.
Let $|A|_k=\sup_{|i-j|=k}|a_{ij}|$ and $|A|_{diag}=\sum_k k|A|_k$.
Then $||[A,P_m]||_1<|A|_{diag}$ for each $m$.
We call $A$ {\it diagonally bounded} if $|A|_{diag}<\infty$. We have
\begin{proposition}
Any finite family of diagonally bounded operators is $1$-semidiagonal.
\end{proposition}

\begin{cor}
Let ${\mathcal A}$ be the algebra of operators on $L_2({\mathbb T})$ generated
by shifts $u(t)\mapsto u(t-\theta)$ and multiplication operators
$M_f$, $f\in C^2({\mathbb T})$. Then any finite family of elements
of ${\mathcal A}$ is $1$-semidiagonal.
\end{cor}
\begin{proof}
It suffices to show that any shift operator and any multiplication operator
$M_f$, $f\in C^2({\mathbb T})$ are diagonally bounded for the standard basis
$e_n=e^{int}$, $n\in{\mathbb N}$. Shifts are diagonally bounded because their
matrices are diagonal. If $f=\sum_na_ne^{int}\in C^2({\mathbb T})$ then
$\sum_n n|a_n|<\infty$ and
$|M_f|_k=max\{|a_k|,|a_{-k}|\}$. Hence $|M_f|_{diag}<\infty$.
\end{proof}
In particular, all Bishop's operators $u(t)\mapsto e^{it}u(t-\theta)$
are $1$-semidiagonal. This was  established by Voiculescu in 
\cite{voi2}.

For a family, ${\mathbb A}=\{A_k\}_{k\in K}$,  of normal 
operators with $\sum_{k\in K}||A_i||^2<\infty$, 
the Hausdorff dimension, \text{dim} of its spectrum 
$\sigma({\mathbb A})\subset l_2(K)$ is appeared
to be important in our study.
We say that the (``essential'') dimension, 
$\text{ess-dim}$, of ${\mathbb A}$ does not exceed $r>0$ if there is
a subset $D$ of $\sigma({\mathbb A})$ such that 
$E_{\mathbb A}(\sigma({\mathbb A})\setminus D)=0$ and $\text{dim}(D)\leq r$
(meaning that there exists $C>0$ such that for $\epsilon>0$ there is
a covering ${\mathcal B}=\{\beta_j\}$ of $D$
by pairwise disjoint Borel sets 
with  $\text{diam}\beta_j<\epsilon$ and 
$|{\mathcal B}|_r:=(\sum_j(\text{diam}\beta_j)^r)^{1/r}\leq C$).
In particular, if $K$ is finite  and all $A_k$ are Lipschitz functions of one 
Hermitian (normal) operator then $\text{ess-dim}({\mathbb A})\leq 1$ 
(respectively $2$). If ${\mathbb A}$ is diagonal then 
$\text{ess-dim}({\mathbb A})=0$.

\begin{proposition}\label{pr53}
If ${\mathbb A}=\{A_k\}_{k\in K}$ is a commutative family of normal operators 
of finite 
multiplicity such that $\text{ess-dim}({\mathbb A})\leq 2$, 
then ${\mathbb A}$ is $2$-semidiagonal.
If it is a finite family of commuting normal operators of finite
multiplicity such that $\text{ess-dim}({\mathbb A})\leq p$, $p<2$, 
then ${\mathbb A}$ is $p$-semidiagonal.
\end{proposition}
\begin{proof}
Suppose first that ${\mathbb A}$ has a cyclic vector. Then  all $A_k$ can be 
realized on $L_2(\sigma({\mathbb A}),\mu)$ as multiplication operators by 
the coordinate functions. Without loss of generality we may assume that
$\text{dim}(\sigma({\mathbb A}))\leq p$. Given a family 
${\mathcal B}=\{\beta_j\}_{j=1}^N$ of
pairwise disjoint Borel subsets of $\sigma({\mathbb A})$ we denote by 
$P_{\mathcal B}$ the projection onto the subspace generated by the
characteristic functions $\chi_j$ of the subsets $\beta_j$.
 Then $$\sum _{k\in K}||[P_{\mathcal B},A_k]||_p^2\leq D|{\mathcal B}|_p^2,$$
where $D$ is a constant.
In fact, let $e_j=\chi_j/||\chi_j||$ and  $\lambda\in\beta_j$. In what follows
we assume $K$ to be finite if $p\ne 2$.  
\begin{eqnarray*}
\displaystyle\sum_{k\in K}||(1-P_{\mathcal B})A_kP_{\mathcal B}e_j||^p\sum_{k\in K}||(1-P_{\mathcal B})(A_k-\lambda_k)P_{\mathcal B}e_j||^p\leq\\
\displaystyle\leq\sum_{k\in K}||(A_k-\lambda_k)P_{\mathcal B}e_j||^p\leq
C(\sum_{k\in K}||(A_k-\lambda_k)P_{\mathcal B}e_j||^2)^{p/2}\leq\\
\displaystyle\leq
C(\sum_{k\in K}|((A_k-\lambda_k)^*(A_k-\lambda_k)E_{\mathbb A}(\beta_j)e_j,e_j)|
)^{p/2}\leq\\
\displaystyle\leq C(\sum_{k\in K}(A_k-\lambda_k)^*(A_k-\lambda_k)E_{\mathbb A}(\beta_j)e_j,e_j))^{p/2}\leq\\
\displaystyle \leq C||\sum_{k\in K}(A_k-\lambda_k)^*(A_k-\lambda_k)E_{\mathbb A}(\beta_j)||^{p/2}
\leq\\
\displaystyle\leq C(\sup_{\lambda'\in\beta_j}\sum_{k\in K}|\lambda'_k-\lambda_k|^2)^{p/2}\leq C(\text{diam}
\beta_j)^p.
\end{eqnarray*}
Here we use the fact that the norms $||\alpha||_p=(\sum_{i=1}^n|\alpha_i|^p)^{1/p}$,
$\alpha=(\alpha_i)\in{\mathbb C}^n$, $1\leq p<\infty$ are equivalent on ${\mathbb C}^n$, $C$
is a  corresponding constant.
Since $\{e_j\}_{j=1}^N$ is a basis of the subspace $P_{\mathcal B}H$ and 
$p\leq 2$ we have
\begin{eqnarray*}
&\displaystyle \sum_{k\in K}||(1-P_{\mathcal B})A_kP_{\mathcal B}||_p^2\leq 
\sum_{k\in K}(\sum_{j=1}^N
||(1-P_{\mathcal B})A_kP_{\mathcal B}e_j||^p)^{2/p}\leq\\ 
&\displaystyle\leq M(\sum_{k\in K}\sum_{j=1}^N
||(1-P_{\mathcal B})A_kP_{\mathcal B}e_j||^p)^{2/p}\leq
MC^{2/p}(\sum_{j=1}^N (\text{diam}\beta_j)^p)^{2/p}=MC^{2/p}|{\mathcal B}|_p^2
\end{eqnarray*}
for some constant $M$ coming from the equivalence of the norms on ${\mathbb C}^n$.

Similarly, $\displaystyle\sum_{k\in K}||(1-P_{\mathcal B})A_k^*P_{\mathcal B}||_p^2\leq
L|{\mathcal B}|_p^2$ ($L=MC^{2/p}$), and therefore  
$\displaystyle\sum_{k\in K}||P_{\mathcal B}A_k(1-P_{\mathcal B})
||_2^2\leq L|{\mathcal B}|_p^2$ so that 
$$\sum _{k\in K}||[P_{\mathcal B},A_k]||_p^2=\sum _{k\in K}||(1-P_{\mathcal B})A_kP_{\mathcal B}||_p^2+\sum_{k\in K}||P_{\mathcal B}A_k(1-P_{\mathcal B})
||_2^2\leq 2L|{\mathcal B}|_p^2.$$
Let now ${\mathcal B}^{(i)}=\{\beta_j^{(i)}\}_{j=1}^{\infty}$ be a sequence of the 
partitions of $\sigma({\mathbb A})$ such that $\sup_j \text{diam}\beta_j^{(i)}\to 0$ and
$|{\mathcal B^{(i)}}|_p\to m_p(\sigma({\mathbb A}))^{1/p}$,
where $m_p$ is the Hausdorff measure on $l_2(K)$. Then $P_{{\mathcal B}^{(i)}}
\to I$ strongly. Since $P_{{\mathcal B}^{(i)}}=s.\lim P_{{\mathcal B}_n^{(i)}}$, where
${\mathcal B}_n^{(i)}=\{\beta_j^{(i)}\}_{j=1}^{n}$, one can find a subsequence of
(finite-dimensional) projections $P_{{\mathcal B}_{n(i)}^{(i)}}$ which converges 
strongly to $I$. 
Finally,
$$\displaystyle\underline{\lim}\sum_{k\in K}||[P_{{\mathcal B}^{(i)}_{n(i)}},A_k]||_p^2\leq
 2L\cdot\underline{\lim}|{\mathcal B}^{(i)}_{n(i)}|_p^2\leq 2L\cdot\lim |{\mathcal B}^{(i)}|_p^2=2L(m_p(\sigma({\mathbb A})))^{2/p}<\infty,$$
implying $p$-semidiagonality of ${\mathbb A}$.

Generally, we decompose the Hilbert space $H$ into a direct sum of subspaces
$H=\oplus_j H_j$, where each $H_j$ is invariant with respect to ${\mathbb A}$ and
${\mathbb A}|_{H_{j}}$ has a cyclic vector. If ${\mathbb A}$ has a finite
multiplicity, we have a finite number of subspaces $H_j$ and the statement easily
follows from what we have already proved.
\end{proof}
%\begin{remark}\rm
%The proposition can be proved for $p<2$ and  infinite families of commuting 
%normal operators ${\mathbb A}=\{A_k\}_{k\in K}$ with finite multiplicity. But in this case
%one has to introduce another definition of the dimension of the family ${\mathbb A}$. 
%\end{remark}

Recall that an operator $A$ is almost normal if $[A^*,A]\in\S_1$. The following result 
was established by Voiculescu \cite[Corollary~2]{voi}.
\begin{proposition}\label{voic}
Any almost normal operator of finite multiplicity is $2$-semidiagonal.
\end{proposition}

In what follows we apply the obtained results to  various problems 
on multiplication operators.

\section{ Application related to the traces of commutators}
%We turn now to deriving corollaries to the above results.

In \cite{weiss} Weiss proved  that if $A$ is a normal operator, $X\in\S_2$ and
$[A,X]\in \S_1$, then $\text{tr}([A,X])=0$. The 
following proposition extends this in several directions.
\begin{proposition}\label{aiipr5}
Let $p\in (1,\infty]$. If $\{A_k\}_{k=1}^n$ is $p/(p-1)$-semidiagonal,
$X_k\in\S_p$ and $\sum_{k=1}^n[A_k,X_k]\in\S_1$ then
$$\tr(\sum_{k=1}^n[A_k,X_k])=0.$$
\end{proposition}
\begin{proof}
Let $T_k:\S_p\to \S_p$, $T_k(X)=[A_k,X]$, and $S_k=T_k|_{\S_1}$.
By Proposition~\ref{aiipr4}(iii) and Remark~\ref{rem},
there exists a common AII for  the intertwinings $(\Phi,S_k,T_k)$.
By assumption, $\sum_{k\in K}[A_k,X_k]=R=\Phi(R)$, 
for some $R\in\S_1$ and
therefore $ R\in\Phi^{-1}(\sum_{k=1}^n \text{Im }T_k)$.
Then Theorem~\ref{aiipr1} gives 
$R\in \overline{(\sum_{k=1}^n\text{Im }S_k)}$ and, since  $\text{Im }S_k\subset 
\ker(\tr)$, we obtain
$\tr(R)=0$.
\end{proof}
\begin{remark}\rm 
The result of Proposition~\ref{aiipr5}  extends to  infinite
family of operators $\{A_k\}_{k\in K}$, $\{X_k\}_{k\in K}$  provided that
that $\sum_{k\in K}||X_k||_p^2<\infty$.
\end{remark}

\begin{cor}
Let  $\{A_k\}_{k=1}^n$ and  $\{B_k\}_{k=1}^n$ be families of operators 
satisfying $$\sum_{k=1}^nB_kA_k=0.$$ 
If $\{A_k\}_{k=1}^n$ is $p/(p-1)$-semidiagonal, $X\in\S_p$ and
$\Delta(X)=\sum_{k=1}^nA_kXB_k\in\S_1$ then
$$\tr(\Delta(X))=0.$$
\end{cor}
\begin{proof}
We have
$$\Delta(X)=\sum_{k=1}^nA_kXB_k-X\sum B_kA_k=\sum_{k=1}^n[A_k,XB_k].$$
Apply now Proposition~\ref{aiipr5} with $X_k=XB_k$.
\end{proof}
\begin{cor}\label{ccor}
If $f_k\in Lip_{1/2}([0,1])$, $1\leq k\leq n$, then no 
functions $F_k\in  L_2([0,1]^2)$ satisfying the condition
\begin{equation}\label{weis}
\sum_{k=1}^n(f_k(x)-f_k(y))F_k(x,y)=1.
\end{equation}
\end{cor}
\begin{proof} Let $X_k$ be the integral operator on $L_2([0,1])$ with the 
kernel $F_k(x,y)$. Clearly, $X_k\in\S_2$. Now 
(\ref{weis}) can be rewritten in the form
$\displaystyle\sum_k[M_{f_k},X_k]=Q$, where $Q$ is the rank-one operator  with
kernel $F(t,s)=1$. By Proposition~\ref{pr53}, the family $\{M_{f_k}\}_{k=1}^n$ is 
$2$-semidiagonal. It remains to apply  Proposition~\ref{aiipr5} to $X_k$, $M_{f_k}$, 
$p=2$, we obtain  $\text{tr }Q=0$. A contradiction.
\end{proof}

One can consider more general classes of functions $F_k$ imposing more
restrictive conditions on $f_k$ (and applying Proposition~\ref{aiipr5} with $p>2$).
For example, for $f_k\in Lip_1[0,1]$, (\ref{weis}) can not hold with
arbitrary function $F_k$ which are the  integral kernels of compact operators.

We do not know if the 
constant $p/(p-1)$ in Proposition~\ref{aiipr5}
is strict for all $p$, but for $p=2$ it is. This follows 
from 
\begin{example}\rm
Let  ${\mathbb D}$ be the unit disk and $dA$ the  area measure on
${\mathbb D}$. In terms of polar coordinates, we have 
$dA(z)=rdrd\theta$, $z=re^{i\theta}$.
Let $H=L^2({\mathbb D}, dA(z))$ and let
$$Au(z)=zu(z),\quad Xu(z)=\int\int(\xi-z)^{-1}u(\xi)dA(\xi).$$
Then $[A,X]=Q$, where $Qu(z)=\int\int u(\xi)dA(\xi)$.
Clearly, $\text{rank }Q=1$, $\text{tr }Q>0$, $A$ is normal and cyclic and hence 
$2$-semidiagonal. We next claim that $X\in\cap_{\epsilon>0}\S_{2+\epsilon}$.
In order to prove this we decompose  $X$ into  the sum
$X=X_1+X_2$, where $X_i$ are defined by similar integrals but 
in $X_1$ we integrate by the disk $|\xi|\leq |z|$, 
in $X_2$ by the annulus $|z|\leq |\xi|\leq 1$. 

Actually, both $X_i$ are represented as operator-weighted bilateral shifts. 
Indeed, let us denote, for any $k\in{\mathbb Z}$, by $H_k$ the subspace of 
$L_2({\mathbb D})$, consisting of functions  $u(r,\theta)=f(r)e^{ik\theta}$,
where $(r,\theta)$ are the polar coordinates.
The map $u\mapsto f$ identifies $H_k$ with the space
$L_2([0,1],dm)$, where $dm=2\pi xdx$. For $u\in H_k$ one has
$$X_1u(z)=\sum_{n=0}^{\infty}\int\int_{r\leq |z|}f(r)e^{ik\theta}r^n
\frac{e^{in\theta}}{z^{n+1}}rdrd\theta=\left\{
\begin{array}{ll}2\pi\int_{r\leq |z|}f(r)r^{-k}
z^{k-1}rdr,&k\leq 0,\\
0,& k>0.
\end{array}\right.$$
Writing $z=re^{i\theta}$ we have therefore
$$X_1u(r,\theta)=h(r)e^{i(k-1)\theta},$$
where $h(r)=A_kf(r)$, $A_k$ is the  integral operator on
$L_2([0,1], dm)$ with  kernel $K(r,t)=r^{k-1}t^{-k}\chi_{t<r}(r,t)$.
One can easily compute that $||A_k||_2^2=\pi^2/(|k|+1)$ so that
for $p>2$
$$||X_1||_p^p=\sum_{k=0}^{\infty}||A_k||_p^p\leq\sum_{k=0}^{\infty}||A_k||_2^p
\leq \sum_{k=0}^{\infty}\pi^p/(|k|+1)^{p/2}<\infty,$$
i.e., $K_1\in\S_p$ for any $p>2$.
Similar arguments shows that $K_2\in\S_p$, $p>2$, verifying the statement.

\end{example}
The above construction answers  a question of Weiss \cite{weiss2} (``does a nuclear
commutator of a compact and a normal operators have zero trace?'').
It was first published in \cite{zap} with the reference to \cite{birman} for the proof 
of the inclusion $X\in\cap_{\epsilon>0}\S_{2+\epsilon}$. We included the proof 
because the 
reference was a mistake - \cite{birman} does not contain this fact.
Using the arguments of [W1] we deduce from the above example an answer to Question 2 in [W3]:

\begin{cor} There is a normal operator $A$ and a compact operator $X$ such that 
$[A,X] \in \S_1$, $[A*,X] \notin \S_1$.
\end{cor}

Weiss \cite{o} asks also  if (\ref{weis}) can be satisfied with 
$n=2$, and $F_k\in L_2([0,1]^2)$ if $f_k$ are only supposed  to be continuous. 
The answer is positive as the following result shows.  

\begin{proposition}
There are functions $f_1$, $f_2\in Lip_{1/3}[0,1]$ and $F_1$, 
$F_2\in L^2([0,1]^2)$ such that 
\begin{equation}\label{lipp}
(f_1(x)-f_2(y))F_1(x,y)+(f_2(x)-f_2(y))F_2(x,y)=1
\end{equation}
\end{proposition}
\begin{proof}
Let $\Pi=[0,1]^3\subset{\mathbb R}^3$, $\pi_i$ be the coordinate functions on 
$\Pi$ ($1\leq i\leq 3$). Set 
\begin{eqnarray*}
&\varphi_1(\lambda)=\pi_1(\lambda),\quad
\varphi_2(\lambda)=\pi_2(\lambda)+i\pi_3(\lambda)\\
&\displaystyle\Phi_1(\lambda,\mu)=\frac{\pi_1(\lambda)-\pi_1(\mu)}{|\lambda-\mu|^2},
\quad
\displaystyle\Phi_2(\lambda,\mu)=\frac{\pi_2(\lambda)-i\pi_3(\lambda)-
\pi_2(\mu)+i\pi_3(\mu)}{|\lambda-\mu|^2}.
\end{eqnarray*}
Then $$\sum_{i=1}^2(\varphi_i(\lambda)-\varphi_i(\mu))\Phi_i(\lambda,\mu)=1.$$
An easy calculation shows that $\Phi_i\in L^2([0,1]^6)$.

Let $\gamma$ be the standard Peano curve $[0,1]\to\Pi$. It can be checked
that $\gamma\in Lip_{1/3}([0,1],{\mathbb R}^3)$ and $\gamma$ preserves the 
Lebesgue measure: $m_1(\gamma^{-1}(E))=m_3(E)$, where $m_n$ is the Lebesgue 
measure on ${\mathbb R}^n$. So, if
we set $f_i(x)=\varphi_i(\gamma(x))$,
$F_i(x,y)=\Phi_i(\gamma(x),\gamma(y))$ then
$F_i\in L_2([0,1]^2,m_2)$, $f_i\in Lip_{1/3}[0,1]$ and (\ref{lipp}) holds.
\end{proof}
Note that if $f_j$ in (\ref{lipp}) are supposed to be real-valued 
then (\ref{lipp}) fails. Indeed, in this case $M_{f_1}$, $M_{f_2}$ is
a pair of commuting self-adjoint operators, hence is $2$-semidiagonal by
Proposition~\ref{pr53}. So it suffices to apply Proposition \ref{aiipr5}. 

\section{Non-commutative version of Fuglede theorem; extensions 
of Fuglede-Weiss theorem}

The well-known problem of the existence of an operator $A$ for which the image
of the derivation $X\mapsto[A,X]$ has non-trivial intersection with
the commutant of the operator $A^*$ can be formulated
for general multiplication operators as the problem
of the validity of the equality
\begin{equation}\label{kerim}
\ker\tilde\Delta\cap \text{Im}\Delta=0
\end{equation}
or of the equivalent equality
\begin{equation}\label{ker}
\ker\tilde\Delta\Delta=\ker\Delta.
\end{equation}

It should be noted that (\ref{ker}) seems to be  the ``right'' form of the 
``Fuglede Theorem for non-normal operators''. Indeed, while the Fuglede
theorem can be considered as the analogs for normal derivations
of the fact that $\ker A=\ker A^*$, for a normal operator $A$, the 
equality  (\ref{ker}) is an analog of the equality
$\ker A^*A=\ker A$, for arbitrary operators. Clearly if
(\ref{ker}) holds and $\Delta$ commutes with $\tilde\Delta$  
 then (\ref{ker}) immediately gives 
$\ker\Delta=\ker\tilde\Delta\Delta=\ker\tilde\Delta$.
Since as we know there are multiplication operators with commuting normal coefficients
for which $\ker\Delta \neq \ker\tilde\Delta$, the equality (\ref{ker}) fails in general. 
Nevertheless the following result shows that it holds for a broad class 
of multiplication operators. 
\begin{theorem}\label{aiicor6}
If  $\{A_k\}_{k\in K}$ is $1$-semidiagonal then (\ref{ker})
is valid. 
\end{theorem}
\begin{proof}
By Corollary~\ref{aiicor4} there exists a $*$-AIS for 
$(\Phi_2,\Delta_{\S_2},\Delta)$ and $(\Phi,\tilde\Delta_{\S_2},\tilde\Delta)$.
The statement now follows from Corollary~\ref{aiicor3} applied 
to $\H=\S_2$, $\Y=B(H)$,
$T_1=\Delta$, $T_2=\tilde\Delta$, $S=\Delta_{\S_2}$.
\end{proof}

\begin{proposition}\label{aiicor7}
If  $\{A_k\}_{k\in K}$ is $1$-semidiagonal and 
$\tilde\Delta\Delta(X)\in\S_1$, then 
$\Delta(X)\in\S_2$.
\end{proposition}
\begin{proof}
By Corollary~\ref{aiicor4},  there exists a $*$-AIS for the pair 
$((\Phi_2,\Delta_{\S_2},\Delta)$,
$(\Phi_2,\tilde\Delta_{\S_2},\tilde\Delta))$.
Moreover,$$\Delta(X)\in\tilde\Delta^{-1}(\S_1)\cap
\text{Im}\Delta=\tilde\Delta^{-1}(\Phi_2\Phi_2^*((B(H)^*))\cap\text{Im}\Delta.$$
Applying now Theorem~\ref{aiipr3}, we obtain $\Delta(X)\in\Phi_2(\S_2)=\S_2$.
\end{proof}

%\begin{proposition} 
%$||\Delta(X)||_2=tr(\tilde\Delta(\Delta(X))X^*)$ for $X?$.
%\end{proposition}
Let us write $||X||_2 = \infty $ if $X \notin \S_2$. In this notation
the famous Fuglede-Weiss
theorem \cite{weiss2} states that 
$$||AX-XB||_2= ||A^*X-XB^*||_2$$
for any normal operators $A,B$ and any operator $X$. Our next task is 
to extend this result to hyponormal operators.

Recall that an operator $A\in B(H)$ is said to be  {\it hyponormal} if
$[A^*,A]$ is positive.

\begin{theorem}\label{aiipr}
Let $A\in B(H)$ be a hyponormal operator of finite multiplicity
and let $B\in B(H)$ be such that $B^*$ is hyponormal. Then for each
$X\in B(H)$
$$||AX-XB||_2\geq ||A^*X-XB^*||_2.$$
\end{theorem} 
\begin{proof}
Let first $X\in\S_2$. Then
\begin{eqnarray*}
&||AX-XB||_2^2=\tr((X^*A^*-B^*X^*)(AX-XB))=\\
&=\tr (X^*A^*AX-AXB^*X^*-XBX^*A^*+
XBB^*X^*)\geq \\
&\geq \tr(X^*AA^*X-AXB^*X^*-XBX^*A^*+XB^*BX^*)=||A^*X-XB^*||_2^2.
\end{eqnarray*}
The inequality for arbitrary $X$  follows now from
Theorem~\ref{aiipr2} applied to ${\mathfrak X}=\S_2$, ${\mathfrak Y}=B(H)$,
$T_1:X\mapsto AX-XB$, $T_2:X\mapsto A^*X-XB^*$ and
$S_i=T_i|_{\S_2}$. The only thing we need to show is the existence
of  AII for
$(\Phi_2,S_1,T_1)$ and, by   
Proposition~\ref{aiipr4}, it would be  sufficient to prove 
that $A$ is $2$-semidiagonal. But beacuse $A$ is hyponormal
and has finite multiplicity,  $A$ is almost normal, i.e. 
$[A,A^*]\in\S_1$ (see, e.g. \cite{clancey}), and therefore  $2$-semidiagonal by
Proposition~\ref{voic}. 
\end{proof}

Now we  consider an extensions of the Fuglede-Weiss theorem to 
general multiplication operators 
with normal coefficients.

Let ${\mathbb A}=\{A_k\}_{k\in K}$ and ${\mathbb B}=\{B_k\}_{k\in K}$,  
be two separately commutating families of normal operators satisfying (\ref{red}).

\begin{proposition}\label{aiicor8}
Suppose that for each $j\in J$ we are given bounded functions $f_j$, $u_j$ on
$\sigma({\mathbb A})$ and $g_j$, $v_j$ on $\sigma({\mathbb B})$ such that 
$$|\sum_{j\in J} f_j(x)g_j(y)|\leq|\sum_{j\in J} u_j(x)v_j(y)|.$$
If $\text{ess-dim }({\mathbb A})\leq 2$, all 
$f_j\in Lip_1(\sigma({\mathbb A}))$ and 
$\sum_{j\in J}||f_j||_{Lip_1}^2<\infty$, then
\begin{equation}\label{ineq}
||\sum_{j\in J}g_j({\mathbb B})Xf_j({\mathbb A})||_2\leq ||\sum_{j\in J}
v_j({\mathbb B})Xu_j({\mathbb A})||_2
\end{equation}
for each $X\in B(H)$.
\end{proposition}
\begin{proof}
Let $\Delta_1(X)=\sum_{j\in J}g_j({\mathbb B})Xf_j({\mathbb A})$ and
$\Delta_2(X)=\sum_{j\in J}v_j({\mathbb B})Xu_j({\mathbb A})$.
The assumption on $\sigma({\mathbb A})$ and the functions $f_j$, 
$j\in J$,
implies $\text{ess-dim}(\{f_j({\mathbb A})\}_{j\in J})\leq 2$.
By Proposition~\ref{pr53}, $\{f_j({\mathbb A})\}_{j\in J}$ is $2$-semidiagonal 
and hence there exists an AII for
$(\Phi_2,(\Delta_1)_{\S_2},\Delta_1)$.
Using the same arguments as in the proof of Proposition~\ref{aiipr}, we see 
that it is enough to show
the inequality (\ref{ineq}) for $X\in\S_2$.
For this we consider the families ${\mathbb A}$, ${\mathbb B}$ concretely
represented as $(A_kf)(x)=a_k(x)f(x)$ on $L_2(T,\mu)$ and
$(B_kg)(x)=b_k(x)g(x)$ on $L_2(S,\nu)$ ($a_k\in L_{\infty}(T,\mu)$, 
$b_k\in L_{\infty}(S,\nu)$). Then $X\in\S_2(L_2(T,\mu ),L_2(S,\nu))$
is an integral operator with a kernel 
$K(x,y)\in L_2(T\times S,\mu\times\nu)$ and  so are the operators 
$\Delta_1(X)$ and $\Delta_2(X)$ with kernels 
$$K_1(x,y)=K(x,y)\sum_{j\in J}f_j(a(x))g_j(b(y)) \text{ and }
K_2(x,y)=K(x,y)\sum_{j\in J}u_j(a(x))v_j(b(y))$$ respectively, where
$a(x)=(a_k(x))_{k\in K}$, $b(x)=(b_k(x))_{k\in K}$.
Since $a(x)\in\sigma({\mathbb A})$, $b(y)\in\sigma({\mathbb B})$ for almost all
$(x,y)\in T\times S$, we have 
$$|K_1(x,y)|\leq |K_2(x,y)|\text{ a.e.},$$
whence
\begin{eqnarray*}
||\Delta_1(X)||_2^2\leq ||\Delta_2(X)||_2^2.
\end{eqnarray*}
\end{proof}

 We mention two special cases of this result; they extend to
"long" multiplication operators the Fuglede-Weiss and Fuglede
theorems respectively. 

%\begin{cor}\label{aiicor9}\cite{}
%Suppose that $A$ and $B$ are normal operators. Then
%$$||f(A)X-Xf(B)||_2\leq ||f||_{Lip_1}\cdot||AX-XB||_2$$
%for any $X\in B(H)$.
%\end{cor}
\begin{cor}\label{aiicor10}
If $\text{ess-dim }({\mathbb A})\leq 2$, then
\begin{equation}\label{wei}
||\sum_{k\in K}B_kXA_k||_2=||\sum_{k\in K}B_k^*XA_k^*||_2
\end{equation}
for any $X\in B(H)$.
\end{cor}

\begin{cor}
Let $A$ be a normal operator,
$f_k\in Lip_1\sigma(A)$,
such that $\displaystyle\sum_{k\in K}||f_k||_{Lip_1}^2<\infty$, and let
$A_k=f_k(A)$. Then the equations
$\displaystyle\sum_{k\in K}A_kXB_k=0$ and
$\displaystyle\sum_{k\in K}A_k^*XB_k^*=0$ are equivalent in $B(H)$.
\end{cor}
\begin{proof}
By the assumptions $\text{ess-dim }\{A_k\}_{k\in K}\leq 2$. 
The statement now trivially follows from Corollary~\ref{aiicor10}.
\end{proof}

It would be desirable to have a ``qualitative'' version of non-commutative 
Fuglede Theorem which would imply simultaneously the Fuglede-Weiss theorem for
a sufficiently general class of multiplication operators. The following
result is a first step in this direction.

\begin{proposition}
If the coefficient family $\{A_k\}_{k\in K}$ of $\Delta$ is $1$-semidiagonal 
then
\begin{equation}\label{quality}
||\Delta(X)||_2^2=tr(X^*\tilde\Delta\Delta(X))
\end{equation}
for any compact operator $X$ such that $\tilde\Delta\Delta(X)\in\S_1$.
\end{proposition}
\begin{proof}
Follows from Theorem~\ref{aiipr3}$(ii)$ and Corollary~\ref{aiicor4}$(i)$.
\end{proof}

Clearly the proposition shows that on $\S_{\infty}$
$\ker\tilde\Delta\Delta=\ker\Delta$ (a restrictive form of 
Theorem~\ref{aiicor6}). On the other hand if $\Delta$ has normal coefficients 
and, for some $X\in\S_{\infty}$, $\tilde\Delta\Delta(X)\in\S_1$ then
$\Delta\tilde\Delta(X)\in\S_1$ and from (\ref{quality}) we get
$||\tilde\Delta(X)||_2^2=||\Delta(X)||_2^2$ - a special case of 
Corollary~\ref{aiicor10}.

Clearly, the ``Fuglede theorem'' for arbitrary $\Delta$ holds in $\S_2$ and
therefore in $\S_p$, $p\leq 2$. Now we will show that for its validity in
$\S_p$, $p>2$, 
some restrictions on the coefficient families
are necessary. 

\begin{proposition}\label{nonsym}
For any $p > 2$ there is a multiplication operator $\Delta$
with commuting normal coefficients satisfying (10) such that
the equations $\Delta(X) = 0$ and $\tilde\Delta(X) = 0$ are
non-equivalent in $\S_p$.
\end{proposition}
\begin{proof}
Using the arguments in \cite[7.8.3-7.8.6]{rudin} one could show that there 
are sequences $c \in l^p(\mathbb Z)$ and $d \in l^1(\mathbb Z)$ such that
$c* d = 0$ and $c*\hat{d} \neq 0$ where $\hat{d}_n = \bar d_{-n}$.
Let $d_n = a_nb_n$ with $|a_n| = |b_n|$, for each $n \in {\mathbb Z}$.
We denote by $U$ the bilateral shift acting on the space $l^2({\mathbb Z})$
($Ue_n = e_{n+1}$, where $(e_n)$ is the standard basis) and set
$A_n = a_nU^n$, $B_n = b_nU^{-n}$, $X = diag(c)$ which means that 
$Xe_n = c_ne_n$. Clearly $\{A_n\}_{n \in \mathbb Z}$,  
$\{B_n\}_{n \in \mathbb Z}$ are commuting families of normal operators 
satisfying (10) and $X \in \S_p$.
An easy calculation shows that $$\sum_{n \in {\mathbb Z}}A_nXB_n = diag(d*c) = 0,\quad 
\sum_{n \in {\mathbb Z}}A_n^*XB_n^* = (diag(d*\hat{c}))^* \neq 0.$$
\end{proof}

A very interesting result of Weiss \cite{weiss1} states that if the 
coefficient families ${\mathbb A}$,
${\mathbb B}$ are (normal and) finite then equality (\ref{wei}) is valid for any $X$ such 
that $\Delta(X)$ and $\tilde\Delta(X)$ belong to $\S_2$. We will finish this 
section by showing how this result may be obtained by using the technique of 
intertwinings. The intermediate steps of the proof are of their own interest and will be
used in the next section.
\begin{lemma}\label{lemma1}
If $\text{card }K<\infty$ then $\S_1\cap\ker\Delta_{\S_2}$ is 
$\S_2$-dense in $\ker\Delta_{\S_2}$.
\end{lemma}
\begin{proof}
Using the spectral theorem we represent our operators $A_k$, $B_k$ as
$(A_kf)(x)=a_k(x)f(x)$ on $H_1=L_2(X,\mu)$ and
$(B_kg)(y)=b_k(y)g(y)$ on $H_2=L_2(Y,\nu)$. Set 
$$\Phi(x,y)=\sum_{k\in K}f_k(x)g_k(y)\text{ and }E=\{(x,y)\in X\times Y\mid \Phi(x,y)=0\}.$$

The space $\S_2(H_1,H_2)$ is naturally identified with 
$L_2(X\times Y,m)$, where $m=\mu\times\nu$. Moreover, one can easily see that
an operator $T$  belongs to  $\ker \Delta_{\S_2}$  iff the corresponding kernel
$K(x,y)\in L_2(X\times Y,m)$ satisfies 
$$K(x,y)\Phi(x,y)=0\quad \text{$m$-almost everywhere.}$$
This means that the space $\ker\Delta_{\S_2}$ is isomorphic to the space 
of all functions $K\in L_2(X\times Y,m)$ vanishing outside  $E$ $m$-almost 
everywhere, i.e. to  the space $L_2(E,m)$.
Since $\text{card K}<\infty$, we have, by Corollary~\ref{scatcor}, that
$E$ is $\tau$-pseudo-open, i.e. it is a union of
a countable number of rectangles $A_i\times B_i$, $A_i\subset X$, 
$B_i\subset Y$ and a $m$-null set. Therefore
$$L_2(E,m)\simeq\oplus_{i=1}^{\infty}L_2(A_i\times B_i,m).$$
It remains to see that for each rectangle $\Pi$ the space $L_2(\Pi,m)$ 
is generated by functions of type $f(x)g(y)$ corresponding to
operators of rank one, the proof is complete.
\end{proof}
Let $\overline{\Delta(B(H))}^{w^*}$ denote the weak$^*$-closure
of $\Delta(B(H))$.
\begin{cor}\label{inters}
If $\text{card }K<\infty$ then 
\begin{equation}\label{formula}
\overline{\Delta(B(H))}^{w^*}\cap\S_2\cap\ker\Delta=\{0\}.
\end{equation}
\end{cor}
\begin{proof}
Let $X\in\overline{\Delta(B(H))}^{w^*}$: $X=\lim_t\Delta(Z_t)$.
For any $Y\in\S_1\cap\ker\Delta$ we have $$\tr(Y^*X)=\lim_t\tr(Z_t(\tilde\Delta(Y))^*).$$
Hence if $X\in\S_2$ then by the previous lemma $\tr(Y^*X)=0$
for any $Y\in\S_2\cap\ker\Delta$. So if also $X\in\ker\Delta$ then 
$\tr(X^*X)=0$ and $X=0$. 
\end{proof}
\begin{cor}\cite{weiss1}. 
If $\text{card }K<\infty$, $\Delta(X),\tilde\Delta(X)\in\S_2$ then
(\ref{wei}) holds.
\end{cor}
\begin{proof}
The equality (\ref{formula}) implies immediately that the 
$\S_2$-closures of $(\text{Im }\Delta)\cap\S_2$ and the $\S_2$-closure of
$(\text{Im }\tilde\Delta)\cap\S_2$  have  trivial intersections with
$\ker\Delta$. Now the result follows directly from Proposition~\ref{intertw}.
\end{proof}

\section{Linear operator equations with normal coefficients.}
The purpose of this section is the study of the thin spectral structure
of the multiplication operators with commuting normal coefficients  satisfying
(\ref{red}). The results will be applied in the next section to
the individual synthesis in Varopoulos
algebras, convolution equations and partial differential equations with constant coefficients.

%$$\sum_{k\in K}B_kXA_k=0,$$
%where ${\mathbb A}=\{A_k\}_{k\in K}$ and ${\mathbb B}=\{B_k\}_{k\in K}$ are commutative families of 
%normal operators acting on $H$ such that
%$\sum_{k\in K}||A_k||^2<\infty$ and $\sum_{k\in K} ||B_k||^2<\infty$.

\begin{lemma}\label{lm}
Let $(\Omega,\mu)$ be a space with finite measure. Assume that $\Omega$ is metrizable 
such that $m_p(\Omega)<\infty$, where $m_p$ is the Hausdorff measure corresponding 
to the measure function $h(r)=r^p$. If $T:L_2(\Omega,\mu)\to H$ is such that
$$||TP(\alpha)||^2\leq C(\text{diam }\alpha)^p$$
for any $\alpha\subset\Omega$, where $P(\alpha)$ is the multiplication operator
by the characteristic function  of $\alpha$, then $T\in\S_2(L_2(\Omega,\mu),H)$ and
$$||T||_2^2\leq C m_p(\Omega).$$
\end{lemma}
\begin{proof}
For any covering ${\mathcal E}=(\alpha_1,\ldots,\alpha_n)$ of $\Omega$,
set $|{\mathcal E}|=\sum_{k=1}^n(\text{diam }\alpha_k)^p$.
Let $e_k=\chi_k/||\chi_k||$, where $\chi_k$ is the characteristic function 
of $\alpha_k$, and let
$Q_{{\mathcal E}}$ be the projection onto the linear span of
$\{e_k\}_{k=1}^n$ in $L_2(\Omega,\mu)$. Then
$$||TQ_{{\mathcal E}}||_{\S_2}=\sum_{k=1}^n||Te_k||^2\leq C|{\mathcal E}|.$$
Taking  a sequence of coverings $\{{\mathcal E}_j\}$ such that
${\mathcal E}_{j+1}$ is a refinement of ${\mathcal E}_j$ and 
$|{\mathcal E}_j|\to m_p(\Omega)$, we obtain
$Q_{{\mathcal E}_j}\to ^s 1$ and
$$\overline{\lim_j}||TQ_{{\mathcal E}_j}||_{\S_2}^2\leq C m_p(\Omega).$$
Therefore, $T\in\S_2$ and $||T||^2_{\S_2}\leq Cm_p(\Omega)$. 
\end{proof}

%Let $S=\{(\lambda,\mu)\in\sigma({\mathbb A})\times\sigma({\mathbb B})\mid
%\sum_{k\in K}\lambda_k\mu_k=0\}$. 

Let $\Delta(X)=\sum_{k\in K}B_kXA_k$, where ${\mathbb A}=\{A_k\}_{k\in K}$, 
${\mathbb B}=\{B_k\}_{k\in K}$ are families of commuting normal operators 
acting
on Hilbert spaces $H_1$ and $H_2$ respectively and satisfying (\ref{strongcoef}).
 Recall that by
${\mathcal E}_{\Delta}(0)$ we denote the space
$$\{T\in B(H_1, H_2)\mid ||\Delta^n(T)||^{1/n}\to 0, n\to\infty\}.$$ 

\begin{lemma}\label{lemma}
Assume that $\text{ess-dim }{\mathbb A}\leq 2n$ and ${\mathbb A}$ has a cyclic 
vector. Then
$$
\Delta^n(X)\in\S_2
$$
for any $X\in{\mathcal E}_{\Delta}(0)$.
\end{lemma}

\begin{proof}
Without loss of generality we can assume that $\text{dim }\sigma({\mathbb A})\leq 2n$. Let $\alpha$ be a closed subset of $\sigma({\mathbb A})$,
$A_k'=A|_{E_{\mathbb A}(\alpha)H_1}$, $\Delta'$ is the multiplication operator with
coefficients $\{A_k'\}_{k\in K}$, $\{B_k\}$ and 
$X'=X|_{E_{\mathbb A}(\alpha)}H_1$. Then it is easy to check that $X'\in{\mathcal E}_{\Delta'}(0)$. Applying  Lemma~\ref{diam} to $\Delta'$, $X'$ we obtain
the equality
$$||\Delta(X)E_{\mathbb A}(\alpha)||\leq
2(\sum_{k\in K}||B_k||^2)^{1/2}||X||(\text{diam }\alpha).$$
Changing repeatedly $X$ by $\Delta(X)E_{\mathbb A}(\alpha)$ we obtain
$$||\Delta^n(X)E_{\mathbb A}(\alpha)||\leq C(\text{diam }\alpha)^{2n},$$
where   $C=(2(\sum_{k\in K}||B_k||^2)^{1/2}||X||)^{2n}$.
 
It remains to note that since $\text{dim }\sigma({\mathbb A})\leq 2n$, we have 
$m_{2n}(\sigma({\mathbb A}))<\infty$, and, since ${\mathbb A}$ has a cyclic 
vector,
$E_{\mathbb A}(\alpha)$  the multiplication by the characteristic functions  of $\alpha$ on
$L_2(\sigma({\mathbb A}),\mu)$, where $\mu$ is the scalar spectral measure of 
${\mathbb A}$.
The statement now follows directly from Lemma~\ref{lm}.
\end{proof}

\begin{theorem}\label{prop}
If $\text{ess-dim }({\mathbb A})\leq 2$, then 
${\mathcal E}_{\Delta}(0)=\ker\Delta$.
\end{theorem}

\begin{proof}
Assume first that ${\mathbb A}$ has a cyclic vector. 
If $X\in {\mathcal E}_{\Delta}(0)$, then, by Lemma~\ref{lemma},
$\Delta(X)\in\Phi_2(\S_2)$. We have 
$$\Delta(X)\in{\mathcal E}_{\Delta}(0)\cap
\S_2={\mathcal E}_{\Delta_{{\mathfrak S}_2}}(0)=\ker\Delta_{\S_2}=\ker(\Delta_{\S_2})^*.$$ The last equality holds because $\Delta_{\S_2}$ is normal.
Therefore, by Corollary~\ref{aiicor1}, 
$\Delta(X)\in\Phi_2(\ker\Delta_{\S_2}^* )\cap\text{Im }\Delta=\{0\}$.

Generally, decompose $H$ into a direct sum of subspaces
 $H=\oplus_{j=1}^{\infty}H_j$, where each $H_j$ is invariant with respect to 
$\{A_k\}_{k\in K}$ and $\{A_k\}_{k\in K}|_{H_j}$ has a cyclic vector.
Then each $X\in B(H)$ can be written as a row operator
$X=(X_1,X_2,\ldots)$, where $X_j=X|_{H_j}$,
and
 $\Delta(X)=(\Delta_1(X_1), \Delta_2(X_2),\ldots)$,
where
$\Delta_j$ is the restriction of $\Delta$ to $B(H_j,H)$.
Now if $X\in{\mathcal E}_{\Delta}(0)$, $X_j\in{\mathcal E}_{\Delta_j}(0)\ker\Delta_j$ and hence $X\in\ker\Delta$. 
\end{proof}

\begin{cor}\label{aiicor12}

If $\text{card }K<\infty$ and $\text{ess-dim }({\mathbb A})\leq 2n$ then
${\mathcal E}_{\Delta}(0)=\ker\Delta^n$
\end{cor}
\begin{proof}
We can assume that $\{A_k\}_{k\in K}$ has a cyclic vector
(the general case reduces to this one as above).
By Lemma~\ref{lemma}, $\Delta^n(X)\in\S_2$ for any 
$X\in {\mathcal E}_{\Delta}(0)$. The arguments similar to one 
in the proof of Theorem~\ref{prop},
 gives 
$$\Delta^n(X)\in\text{Im }\Delta^n\cap\ker\Delta^*_{\S_2}\subset
\text{Im }\Delta\cap\ker\Delta^*_{\S_2}.$$
Therefore it is enough to show that
$\text{Im }\Delta\cap\ker\Delta^*_{\S_2}=\{0\}$.
 But this follows immediately from Corollary \ref{inters}. 
\end{proof}

Let $\X$ be a Banach space and let $T$ be a linear mapping on $\X$. 
The ascent, $\text{asc }T$,  is defined as the least positive integer $n$ such that
$\ker T^n=\ker T^{n+1}$. 
If no such  integers exist   we put
$\text{asc }T=\infty$.
Since $\ker\Delta^k\subset{\mathcal E}_{\Delta}(0)$,
we obtain from Corollary~\ref{aiicor12} an estimate of the ascent
\begin{equation}
\text{ess-dim }({\mathbb A})\leq 2n\Rightarrow 
\text{asc }\Delta\leq n
% \frac{1}{2}(\text{dim }\sigma({\mathbb A})) \quad 
%(\text{asc }\Delta\leq\frac{1}{2}(\text{dim }\sigma({\mathbb A})+1)????)
\end{equation}
or setting $(x]=-[-x]$, the smallest integer $\geq x$,
\begin{equation}\label{ascf}
\text{asc }\Delta\leq (\frac{1}{2}\text{ess-dim }({\mathbb A})]
\end{equation}
This implies a more simple and rough result:
$\text{asc }\Delta \le k = \text{card }K$. Somewhat 
more precise (but of course also rough) estimate is given in the following
result:
\begin{proposition}\label{asc}
$\text{asc }\Delta \le k-1$.
\end{proposition}
\begin{proof}
Clearly,   if $\Delta$ has length $k$ then, by (\ref{ascf}),
 $\text{asc }(\Delta-z I)\leq k$ for any constant $z$.

Assume first that the operators $A_1$, $B_1$ are invertible and denote by 
$R_{A_1}$, $L_{B_1}$  the right and the left multiplication  by
$A_1$ and $B_1$ respectively.
Then
$$\Delta=R_{A_1}L_{B_1}(\Delta'+1),$$
where $\Delta'$  is a multiplication operator of length $\leq k-1$.
Since clearly $\text{asc }\Delta=\text{asc }(\Delta'+1)$, we obtain 
$\text{asc }\Delta\leq k-1$.

To prove the statement in general case consider the hyperplane
$S=\{(z_i)\in{\mathbb C}^k\mid z_1=0\}$ and the set $M$
of all closed subsets $K\subset{\mathbb C}^k$ such that either $K\cap S=\emptyset$
or $K=S$. Let $Q$ be the family of the projections $R_{E_{\mathbb A}(K_1)}
L_{E_{\mathbb B}(K_2)}$, $K_1$, $K_2\in M$. One can easily see
that $Q$ is complete, meaning that, for any $X\in B(H_1,H_2)$,  the closed subspace
generated by $P(X)$, $P\in Q$, contains $X$.

Next we note that 
$$\text{asc }\Delta=\sup\{\text{asc}(\Delta P)\mid P\in Q\}$$
for any complete family $Q$ of projections commuting with $\Delta$ and that
for any such projection $P$ either $\text{asc}(\Delta P)=1$ or $\text{asc}(\Delta P)\text{asc}(\Delta P|_{P\X})$, $\X=B(H_1,H_2)$.
Hence it is enough to show that 
$\text{asc }(\Delta P|_{P\X})\leq n-1$ for any $P\in Q$. But if some of $K_1$ and $K_2$
equals $S$ then for $P=R_{E_{\mathbb A}(K_1)}
L_{E_{\mathbb B}(K_2)}$, $R_{A_1}L_{B_1}P=0$ and $\Delta P|_{P\X}$ has length $\leq n-1$ 
implying  $\text{asc }(\Delta P|_{P\X})\leq n-1$. Otherwise, the restrictions of the 
operators  $R_{A_1}$, $L_{B_1}$ to
$P\X$ are invertible reducing the problem to the case treated in the beginning.
\end{proof}
\begin{remark}\rm Proposition~\ref{asc}, being much less general than (\ref{ascf}) has 
its advantages. For example, it implies immediately the result of Weiss \cite{weiss2}  
on multiplication operators of the length 2.
\end{remark}
\begin{cor}\label{lip}
If $\text{ess-dim }({\mathbb A})\leq 2$, then the solution space of the equation
$$\sum_{k\in K}B_kXA_k=0$$ is reflexive.
\end{cor}
\begin{proof}
It follows from Theorem~\ref{prop}, Proposition~\ref{S} and 
Corollary~\ref{co}.
\end{proof}
In the next section it  
will be shown that this statement can be regarded as an operator version of 
the Beurling-Pollard theorem on synthesis of $Lip_{1/2}$-functions
on the circle.

\section{Individual synthesis in Varopoulos algebras; some applications}
We can  now return to our initial topic and apply 
 Proposition~\ref{indiv} and Theorem~\ref{prop} to obtain a criterium for
synthesis of functions in the Varopoulos algebra $V(X,Y)$.
\begin{theorem}\label{TBP}
Let $F=\sum_{i=1}^{\infty} f_i\otimes g_i\in V(X,Y)$, 
$f(x)=(f_1(x),f_2(x),\ldots)\in l_2$.
If $\text{dim } f(X)\leq 2$, then 
$F$ admits spectral synthesis.
\end{theorem}
\begin{proof}
Take arbitrary borel measures $\mu$, $\nu$ on $X$, $Y$ and let  $\Delta_F$ be the multiplication operator as defined in
Section~\ref{tensor}. Then its left coefficient family
${\mathbb A}=\{A_i\}_{i=1}^{\infty}$ coincides with 
$\{M_{f_i}\}_{i=1}^{\infty}$, 
whence 
$\sigma({\mathbb A})\subset f(X)$, 
$\text{dim }\sigma({\mathbb A})\leq\text{dim }f(X)\leq 2$.
The statement now follows from Theorem~\ref{prop} and Corollary~\ref{co}.
\end{proof}
Note that
if $\text{dim }X=1$ (or $2$)  and $f_i\in Lip_{1/2}(X)$ 
(respectively $f_i\in Lip(X)$) 
with the Lipschitz constants $C_i$ such that 
$\sum C_i <\infty$ the  theorem  says that 
 $F(x,y)=\sum_if_i(x)g_i(y)$ admits spectral synthesis in $V(X,Y)$. 
This shows that 
 Theorem~\ref{TBP} can be considered as a tensor algebra version of the 
Beurling-Pollard theorem and 
Corollary~\ref{lip} as its operator version.
 
\begin{theorem}
Let $F(x,y)=\sum_{i=1}^k f_i(x)g_i(y)\in V(X,Y)$.
Let $m=(\text{dim } f(X)/2]$, where
$f:X\to{\mathbb C}^k$ is the map $x\mapsto(f_1(x),\ldots, f_k(x))$.
Then the sequence of closed ideals $J_j=\overline{F^jV(X,Y)}$ of
$V(X,Y)$ stabilizes on a number $n\leq m$. Moreover, for any
Banach module $M$ over $V(X,Y)$ the sequence of submodules 
$\overline{F^jM}$ stabilizes on a number $n\leq m$.
\end{theorem}
\begin{proof}
Clearly the annihilator $J_j^{\perp}$ of $J_j$ in $V(X,Y)'$
coincides with the space

$$W_j=\{B\in V(X,Y)'\mid F^jB=0\}.$$ 
It suffices to
prove that $W_{m+1}=W_m$. Let $B\in W_{m+1}$.
As in the proof of Proposition~\ref{indiv} there are measures $\mu$, $\nu$
on $X$, $Y$ and an operator $T:L_2(X,\mu)\to L_2(Y,\nu)$ such that
$$\langle u\otimes v,B\rangle=(Tu,v)$$
for any $u\in C(X)$, $v\in C(Y)$. 

If  $\Delta_F$ is the multiplication operator corresponding to $F$
 then for its left coefficient family
${\mathbb A}=\{M_{f_i}\}_{i=1}^k$, we have
$\sigma({\mathbb A})\subset f(X)$, and $\text{dim }\sigma({\mathbb A})\leq\text{dim }f(X)$.
The bimeasure $FB$ corresponds to the 
operator $\Delta_F(T)$, $F^jB$ to $\Delta_F^j(T)$.
Thus $T\in\ker\Delta_{F^{m+1}}=\ker\Delta_{F^{m}}$ by (\ref{ascf}),
and $B\in W_m$.

We actually established that $F^m=\lim_jF^{m+1}G_j$ for
some sequence $\{G_j\}$ in $V(X,Y)$. This implies immediately the equality
$\overline{F^mM}=\overline{F^{m+1}M}$ for any Banach module $M$.
\end{proof}

Next results relate the problematic with harmonic analysis and ordinary
differential equations.
Let ${\mathcal F}$ denote the Fourier transform  
(acting in any of the spaces of
ordinary or generalized functions considered below), $D$ be the space
of compactly supported infinitely differentiable functions  on
${\mathbb R}^n$, $D'$ its dual 
space, ${\mathcal F}L_1({\mathbb R}^n)$ the Fourier algebra, 
$PM({\mathbb R}^n)$ the space dual to ${\mathcal F}L_1({\mathbb R}^n)$
(the space of pseudomeasures), $\varphi*\psi$ the convolution of two functions
 in $D$. The imbedding $PM\subset D'$ permits one to consider the
distribution $p\Phi\in D'$ for any polynomial $p$ in $n$ variables
and  any pseudomeasure $\Phi$.
 
\begin{cor}\label{corhar}
Let $p$ be a polynomial in two variables then for $\Phi\in PM({\mathbb R}^2)$
the inclusion $\supp\Phi\subset p^{-1}(0)$ is equivalent to the condition
$p\Phi=0$. 
\end{cor}
\begin{proof}
Let $\supp\Phi\subset p^{-1}(0)$. Clearly, there exist polynomials
$s_i$, $r_i$, $1\leq i\leq N$, such that $p(x-y)=\sum_{i=1}^Ns_i(x)r_i(y)$,
$x$, $y\in{\mathbb R}^2$. For $u$, $v\in D$ set $a_i(x)=u(x)s_i(x)$ and
$b_i(y)=v(y)r_i(y)$. Since the Fourier transform ${\mathcal F}\Phi$ of a 
pseudomeasure $\Phi$ belongs to $L_{\infty}({\mathbb R}^2)$ we have
a well-defined bounded operator 
$T={\mathcal F}^{-1}M_{{\mathcal F}\Phi}{\mathcal F}$ on $L_2({\mathbb R})$,
here $M_{{\mathcal F}\Phi}$ is the multiplication operator by the function
${\mathcal F}\Phi$.
$T$ is supported in $$E=\{(x,y)\in{\mathbb R}^2\times{\mathbb R}^2\mid
x-y\in \text{supp }\Phi\}.$$
In fact, for $f$, $g\in L_2({\mathbb R}^2)$, one can easily see that
$$(Tf,g)=\langle\Phi,f*\tilde{\bar g}\rangle,$$
where $\tilde\phi(x)=\phi(-x)$. Let $\alpha$, $\beta$ be closed sets such that
$(\alpha\times\beta)\cap E=\emptyset$. Since $E$ is closed, there
exist open sets $\alpha_0\supset\alpha$, $\beta_0\supset\beta$ such that
$\bar\alpha_0\times\bar\beta_0$ does not intersect $E$. For every
pair of functions  
$f$, $g\in L_2({\mathbb R}^2)\cap C({\mathbb R}^2)$, which vanish 
outside the sets $\alpha_0$ and $\beta_0$ respectively, we have 
$\supp(f*\tilde{\bar g})\subset\supp(f)-\supp(g)\subset
\bar\alpha_0-\bar\beta_0$, and since $(\bar\alpha_0-\bar\beta_0)\cap\supp\Phi=\emptyset$,
we have $(Tf,g)=\langle\Phi,f*\tilde{\bar g}\rangle=0$,
implying $M_{\chi_{\beta_0}}TM_{\chi_{\alpha_0}}=0$ and 
$M_{\chi_{\beta}}TM_{\chi_{\alpha}}=0$. By the regularity of the Lebesgue 
measure we have that this is true for any Borel sets $\alpha$, $\beta$.
Clearly 
 $$E\subset\{(x,y)\in{\mathbb R}^2\times{\mathbb R}^2\mid p(x-y)=0\}\subset
\{(x,y)\in{\mathbb R}^2\times{\mathbb R}^2\mid \sum_{i=1}^Na_i(x)b_i(y)=0\}.$$
Since  $a_i$ are smooth functions on ${\mathbb R}^2$, we have that
$\text{ess-dim }(\{M_{a_i}\}) \le 2$.
Applying now Proposition~\ref{ref1} and Theorem~\ref{prop} we conclude that
$\sum M_{b_i}TM_{a_i}=0$. 
A direct computation shows that for $\varphi$, $\psi\in D$
$$(\sum_{i=1}^N M_{b_i}TM_{a_i}\varphi,\psi)=\langle p\Phi,u\varphi *
\widetilde{\bar{v}\psi}\rangle,$$
where $\langle\cdot,\cdot\rangle$  is the  pairing of
the spaces $D$ and $D'$. Therefore,
$\langle p\Phi,u\varphi*
\widetilde{\bar{v}\psi}\rangle  = 0$, for  any $u$, $\varphi$, $v$, 
$\psi\in D$; this shows that 
$p\Phi=0$.

The reverse implication is obvious.
\end{proof}

\begin{cor}\label{diffeq}
The space of all bounded solutions of the equation
\begin{equation}\label{diff}
p(i\frac{\partial}{\partial x_1},i\frac{\partial}{\partial x_2})u=0
\end{equation}
in ${\mathbb R}^2$ is completely determined by the variety of zeros of
the polynomial $p$ in ${\mathbb R}^2$.
\end{cor}

\begin{proof}
The equation (\ref{diff}) is equivalent to
$p\Phi=0$, where $\Phi={\mathcal F}^{-1}u$ is a pseudomeasure.
By Corollary~\ref{corhar}, the space of its solutions  is the set
of all pseudomeasures such that $\supp\Phi\subset p^{-1}(0)$.
\end{proof}

\begin{remark}\rm
(i) The result of Corollary~\ref{diffeq} clearly extends to a wide class of 
infinite order equations (that is (\ref{diff}) with a smooth functions $p$ 
instead of polynomial).

(ii) For polynomials in $n>2$ variables the result obtains the following
form: if $p$, $q$ are polynomials with $p^{-1}(0)\subset q^{-1}(0)$ then any
bounded solution of the equation $$p(i\frac{\partial}{\partial x_1},\ldots,
i\frac{\partial}{\partial x_n})u=0$$
satisfies the equation
$$q(i\frac{\partial}{\partial x_1},\ldots,
i\frac{\partial}{\partial x_n})^mu=0,$$
where $m=(n/2]$. For the proof one should only use Corollary~\ref{aiicor12}
instead of Theorem~\ref{prop}.
\end{remark}

\begin{center}
{\footnotesize \sl Department of Mathematics, Vologda State Liceum of
Mathematical and Natural Sciences,
%}\\
%{\footnotesize \sl
Vologda, 160000, Russia}\\
{\footnotesize \sl
 shulman\_v\symbol{64}yahoo.com}\\
 \vspace{0.5cm}
{\footnotesize \sl Department of Mathematics, Chalmers University of 
Technology,  
SE-412 96 G\"oteborg, Sweden} \\
{\footnotesize \sl
turowska\symbol{64}math.chalmers.se}
\end{center}
\end{document}